\documentclass[11pt]{article}
\usepackage{amssymb}
\usepackage{latexsym}
\usepackage{amsfonts}
\usepackage{amsmath}
\newtheorem{thrm}{Theorem}[section]
\newtheorem{lemma}[thrm]{Lemma}

\newtheorem{remark}[thrm]{Remark}
\newtheorem{example}[thrm]{Example}
\numberwithin{equation}{section}

\usepackage[dvips]{color}

\def\cL{{\cal L} }
\def\P{\mathbb{P} }
\def\E{\mathbb{E} }

\begin{document}
\allowdisplaybreaks
\begin{doublespace}

\title{\Large\bf Central Limit Theorems for Supercritical Branching Markov Processes}
\author{ \bf  Yan-Xia Ren\footnote{The research of this author is supported by NSFC (Grant No.  11271030 and 11128101) and Specialized Research Fund for the Doctoral Program of Higher Education.\hspace{1mm} } \hspace{1mm}\hspace{1mm}
Renming Song\thanks{Research supported in part by a grant from the Simons
Foundation (208236).} \hspace{1mm}\hspace{1mm} and \hspace{1mm}\hspace{1mm}
Rui Zhang\footnote{Supported by the China Scholarship Council.}
\hspace{1mm} }
\date{}
\maketitle

\begin{abstract}
In this paper we establish spatial central limit theorems for a large class of
supercritical branching Markov processes with general spatial-dependent branching mechanisms.
These are generalizations of the spatial central limit theorems proved
in \cite{RP}
for branching OU processes with binary branching mechanisms. Compared with the results of
\cite{RP},
our central limit theorems are more
satisfactory in the sense that the normal random variables in our theorems are non-degenerate.
\end{abstract}


\medskip


\section{Introduction}

In recent years, there have been many papers on law of large numbers type convergence theorems
for branching Markov processes and superprocesses, see, for instance,
\cite{CRW, CS, E09, EHK, EW, LRS2, MR, W} and the references therein. For recent results
on other non-central limit theorem types convergence results for branching Markov processes,
see, for instance, \cite{Ha, HH, LRS, LRS11}
and the references therein.

The focus of this paper is on spatial central limit theorems for branching Markov processes.
For critical branching Markov processes starting from a Poisson random field or an equilibrium distribution,
and subcritical branching Markov processes with immigration,
some functional central limit theorems of the occupation times were established in a series of papers,
see, for instance,
\cite{BGT06a, BGT06b, BGT08, Mi08, Mi09, Mi12}
and reference therein.
However, up to now, no spatial central limit theorems have been established for
$\langle f, X_t\rangle$ of
general supercritical  branching Markov processes starting from general initial configurations.
In \cite{RP}, some spatial central limit theorems were established for
$\langle f, X_t\rangle$ of
supercritical branching OU processes with binary branching mechanism starting from a point mass.
In \cite{Mi}, some spatial central limit theorems were established for
supercritical super-OU processes with binary branching mechanisms
starting from finite and compactly supported measures.
However, the central limit theorems of \cite{RP, Mi}
are not very satisfactory since the limiting normal random variables maybe degenerate. In the recent preprint
\cite{RSZ}, we established spatial central limit theorems for
supercritical super-OU processes with general branching mechanisms
starting from finite and compactly supported measures.
The limiting normal random variables in our central limit theorems are non-degenerate.
For earlier central limit theorems for
supercritical branching processes and supercritical multi-type branching processes, see
\cite{Ath69a, Ath69, Ath71, KS66}.

In this paper, we will extend the arguments of \cite{RP, Mi, RSZ} to establish spatial central limit
theorems for a large class of supercritical branching Markov processes
with general spatial-dependent branching mechanisms.

\subsection{Spatial process}\label{subs:sp}

In this subsection, we spell out our assumptions on the spatial Markov process and then give some examples.

Suppose that $E$ is a locally compact separable metric space and that $\mu$ is a $\sigma$-finite Borel measure
on $E$ with full support.
Suppose that $\partial$ is a separate point not contained
in $E$. $\partial$ will be interpreted as the cemetery point.
We will use $E_{\partial}$ to denote $E\cup\{\partial\}$.
Every function $f$ on $E$ is automatically extended to $E_{\partial}$ by setting $f(\partial)=0$.
 We will assume that $\xi=\{\xi_t,\Pi_x\}$ is a $\mu$-symmetric Hunt process on $E$ and $\zeta:=
 \inf\{t>0: \xi_t=\partial\}$ is the lifetime of $\xi$.
We will use $\{P_t:t\geq 0\}$ to denote the semigroup of $\xi$.  Our standing assumption on $\xi$ is
that there exists a family of
continuous strictly positive symmetric functions $\{p_t(x,y):t>0\}$ on $E\times E$ such that
$$
  P_tf(x)=\int_E p_t(x,y)f(y)\mu(dy).
$$
It is well-known and easy to check that, for $p\geq 1$, $\{P_t:t\ge 0\}$ is a
strongly continuous contraction semigroup on $L^p(E, \mu)$.
In fact, it follows from H\"{o}lder's inequality, Fubini's theorem and symmetry that
$$
  \|P_tf\|_p^p=\int_E \left|\int_E p_t(x,y)f(y)\mu(dy)\right|^p\mu(dx)\leq \int_E\int_E p_t(x,y)|f|^p(y)\mu(dy)\mu(dx)\leq \|f\|_p^p.
$$

Define $\widetilde{a}_t(x):=p_t(x, x)$.
Throughout this paper, we will assume that $\widetilde{a}_t(x)$ satisfies the following two conditions:
\begin{description}
  \item[(a)] For any $t>0$, we have
$$
        \int_E \widetilde{a}_t(x)\mu(dx)<\infty.
$$
   \item[(b)] There exists $t_0>0$ such that $\widetilde{a}_{t_0}(x)\in L^2(E,\mu)$.
\end{description}
It is well-known (see, for instance, \cite[Section 2]{DS}) that
$p_t(x, y)\le (\widetilde{a}_t(x)\widetilde{a}_t(y))^{1/2}$
and that, for each $x\in E$, the function $t\to \widetilde{a}_t(x)$ is a decreasing function.
So condition $(b)$ above is equivalent to
 \begin{description}
   \item[(b$'$)] There exists $t_0>0$ such that for all $t\ge t_0$, $\widetilde{a}_{t}(x)\in L^2(E,\mu)$.
 \end{description}
Now we give some examples of Markov processes satisfying the above assumptions.
The purpose of these examples is to show that the above assumptions are satisfied by
many Markov processes. We will not try to give the most general examples possible.
The first example below contains OU processes as special cases.

\begin{example}[Subordinate OU Process]
{\rm
Let $\sigma,b>0$ be two constants.
Suppose that
$\eta=\{\eta_t: t\ge0\}$
is an Ornstein-Uhlenbeck process (OU process, for short)
on $\mathbb{R}^d$, that is, a diffusion process with infinitesimal generator
$$
  L:=\frac{1}{2}\sigma^2\bigtriangleup-b x\cdot\bigtriangledown.
$$
For any $x\in \mathbb{R}^d$, we use $\Pi_x$ to denote the law of $\xi$ starting from $x$.
It is well known that under $\Pi_x$, $\eta_t\sim\mathcal{N}(xe^{-b t}, \sigma_t^2)$,
where $\sigma_t^2=\sigma^2(1-e^{-2b t})/(2b)$ and $\eta$ has an invariant density
$$
\mu(x)=\left(\frac{b}{\pi\sigma^2}\right)^{d/2}\exp \left(-\frac{b}{\sigma^2}\|x\|^2\right).
$$
Let
$$
p^0_t(x,y):=\left(\frac{1}{2\pi\sigma_t^2}\right)^{d/2}\exp \left(-\frac{\|y-xe^{-bt}\|^2}{2\sigma_t^2}\right).
$$
So
\begin{equation}\label{2.8}
  p^0_t(x,x)=\left(\frac{1}{2\pi\sigma_t^2}\right)^{d/2}
  \exp \left(-\frac{b(1-e^{-bt})}{\sigma^2(1+e^{-bt})}\|x\|^2\right).
\end{equation}
Put $E=\mathbb{R}^d$ and $\mu(dx)=\mu(x)dx$. The density  of $\eta_t$ with respect to $\mu$ is
$$
\overline{p}^0_t(x,y)=p^0_t(x,y)\mu(y)^{-1}=
\left(\frac{1}{1-e^{-2bt}}\right)^{d/2}\exp\left\{-\frac{b}{\sigma^2(e^{2bt}-1)}
\left(\|y\|^2+\|x\|^2-2x\cdot ye^{bt}\right)\right\}.$$
In particular,
$$
   \overline{p}^0_t(x,x)=
  \left(\frac{1}{1-e^{-2bt}}\right)^{d/2}\exp\left\{\frac{2b}{\sigma^2(e^{bt}+1)}\|x\|^2\right\}.
$$

Suppose that $S_t$ is a subordinator, independent of $Y$, with Laplace exponent $\varphi$, that is,
$$
 \E(e^{-\theta S_t})=e^{-t\varphi(\theta)}, \qquad \theta>0.
$$
Suppose that $S$ has a positive drift coefficient $a>0$.
Then $S_t\geq at$, for all $t>0$.

The process $\{\xi_t:t\geq0\}$ defined by $\xi_t:=\eta_{S_t}$ is called a subordinate OU process. In the special
case $S_t\equiv t$, $\xi$ reduces to the OU process $\eta$.
Thus the transition density of $\xi_t$ with respect to $\mu$ is given by
$$
p_t(x,y)=
\E\left(\overline{p}^0_{S_t}(x,y)\right).
$$
So $p_t(x,y)$ is symmetric.
By \eqref{2.8}, we have
$$
\int_E \widetilde{a}_t(x)\mu(dx)=\E \int_E p^0_{S_t}(x,x)dx
=\E(1-e^{-bS_t})^{-d}\leq (1-e^{-abt})^{-d}<\infty.
$$
Chose $t_0>0$ such that $4/(e^{abt_0}+1)<1.$
Then by H\"{o}lder's inequality, we get
\begin{eqnarray*}
  \int_E \widetilde{a}_{t_0}^2(x)\mu(dx) &\leq &
  \E\left(\int_E \overline{p}^0_{S_{t_0}}(x,x)^2\mu(x)\,dx\right).
\end{eqnarray*}
For $t\geq at_0$, we have
\begin{eqnarray*}
  \int_E \overline{p}^0_t(x,x)^2\mu(x)\,dx &=&
  \int_{\mathbb{R}^d}\left(\frac{b}{\pi\sigma^2(1-e^{-2bt})^2}\right)^{d/2}
  \exp\left\{-\left(1-\frac{4}{e^{bt}+1}\right)\frac{b}{\sigma^2}\|x\|^2\right\}\,dx\\
   &=& \left((1-e^{-2bt})^2\left(1-\frac{4}{e^{bt}+1}\right)\right)^{-d/2}\\
   &\leq& \left((1-e^{-2abt_0})^2\left(1-\frac{4}{e^{abt_0}+1}\right)\right)^{-d/2},
\end{eqnarray*}
which implies
$$
\int_E \widetilde{a}_{t_0}^2(x)\mu(dx)\leq
\left((1-e^{-2abt_0})^2\left(1-\frac{4}{e^{abt_0}+1}\right)\right)^{-d/2}<\infty.
$$
Thus the process $\xi$ satisfies all the assumptions in the beginning of this subsection.
}
\end{example}

\begin{example} {\rm
Suppose $a>2$ is a constant. Let $\xi$ be a Markov process on $\mathbb{R}^d$ corresponding to the infinitesimal
generator $\Delta-|x|^a$. Let $p_t(x, y)$ denote the transition density of $\xi$ with respect to
the Lebesgue measure on $\mathbb{R}^d$.
It follows from \cite[Section 4.5]{D} that, for any $t>0$, there exists $c_t>0$ such that
$$
p_t(x, y)\le c_t\exp\left(-\frac{2}{2+a}|x|^{1+a/2}\right)\exp\left(-\frac{2}{2+a}|y|^{1+a/2}\right),
\qquad x, y\in \mathbb{R}^d.
$$
Taking $E=\mathbb{R}^d$ and $\mu$ to be the Lebesgue measure on $\mathbb{R}^d$, using the display above,
one can easily check that all the assumptions at the beginning of this subsection  are satisfied in this case.
}
\end{example}

\begin{example}{\rm
Suppose that $V$ is a nonnegative and locally bounded function on  $\mathbb{R}^d$
such that there exist $R>0$ and $M\ge 1$ such that for all $|x|>R$,
$$
M^{-1}(1+V(x))\le V(y)\le M(1+V(x)), \qquad y\in B(x, 1),
$$
and that
$$
\lim_{|x|\to\infty}\frac{V(x)}{\log|x|}=\infty.
$$
Suppose $\alpha\in (0, 2)$ is a constant.
Let $\xi$ be a Markov process on $\mathbb{R}^d$ corresponding to the infinitesimal
generator $-(-\Delta)^{\alpha/2}-V(x)$. Let $p_t(x, y)$ denote the transition density of $\xi$ with respect to
the Lebesgue measure on $\mathbb{R}^d$.
It follows from \cite[Corollaries 3 and 4]{KK} that, for any $t>0$, there exists $c_t>0$ such that
$$
p_t(x, y)\le c_t\frac1{(1+V(x))(1+|x|)^{d+\alpha}}\frac1{(1+V(y))(1+|y|)^{d+\alpha}},
\qquad x, y\in \mathbb{R}^d.
$$
Taking $E=\mathbb{R}^d$ and $\mu$ to be the Lebesgue measure on $\mathbb{R}^d$, using the display above,
one can easily check that all the assumptions at the beginning of this subsection  are satisfied in this case.
}
\end{example}

\begin{example}{\rm
A nondecreasing function $L:[0, \infty)\to [0, \infty)$ is said to be in the class ${\bf L}$
if $\lim_{t\to\infty}L(t)=\infty$ and there exists $c>1$ such that
$$
L(t+1)\le c(1+L(t)), \qquad t\ge 0.
$$
Suppose that $V$ is a nonnegative function on  $\mathbb{R}^d$
such that
$$
\lim_{|x|\to\infty}\frac{V(x)}{|x|}=\infty
$$
and that there exists a function $L\in {\bf L}$ such that there exists  $C>0$ such that
$$
L(|x|)\le V(x)\le C(1+L(|x|), \qquad x\in \mathbb{R}^d.
$$
Suppose that $m>0$ and $\alpha\in (0, 2)$ are constants.
Let $\xi$ be a Markov process on $\mathbb{R}^d$ corresponding to the infinitesimal
generator $m- (-\Delta +m^{2/\alpha})^{\alpha/2}-V(x)$.
Let $p_t(x, y)$ denote the transition density of $\xi$ with respect to
the Lebesgue measure on $\mathbb{R}^d$.
It follows from \cite[Theorem 1.6]{KuSi} that, for any $t>0$, there exists $c_t>0$ such that
$$
p_t(x, y)\le c_t\frac{\exp(-m^{1/\alpha}|x|)}{(1+V(x))(1+|x|)^{(d+\alpha+1)/2}}
\frac{\exp(-m^{1/\alpha}|y|)}{(1+V(y))(1+|y|)^{(d+\alpha+1)/2}},
\qquad x, y\in \mathbb{R}^d.
$$
Taking $E=\mathbb{R}^d$ and $\mu$ to be the Lebesgue measure on $\mathbb{R}^d$, using the display above,
one can easily check that all the assumptions at the beginning of this subsection  are satisfied in this case.
}
\end{example}

The next example shows that a lot of  important
Markov processes on bounded subsets of $\mathbb{R}^d$ satisfy the above assumptions.

\begin{example}{\rm
Suppose that $E$ is a locally compact separable metric space, $\mu$ is a finite Borel measure
on $E$ with full support and that $\xi=\{\xi_t, \Pi_x\}$ is a $\mu$-symmetric
Hunt process on $E$. Suppose that, for each $t>0$, $\xi_t$ has a continuous, symmetric and strictly
positive density $p_t(x,y)$ with respect to $\mu$. If the semigroup of $\xi$ is ultracontractive,
or equivalently, for
any $t > 0$, there exists constant $c_t > 0$ such that
$$
  p_t(x,y)\le c_t,\quad \mbox{for any }(x,y)\in E\times E.
$$
Then it is trivial to see that, in this case,
all the assumptions at the beginning of this subsection are satisfied.

Some particular cases of this example are as follows:

\begin{description}
\item{(1)}
Suppose that $D$ is a connected open subset of  $\mathbb{R}^d$ with finite Lebesgue measure
and that $\mu$ denotes the Lebesgue measure on $D$. Then the subprocess in $D$ of any diffusion
process in $\mathbb{R}^d$ corresponding to a uniformly elliptic divergence form second order
differential operator satisfies the assumptions of the first paragraph in this example and therefore
all the assumptions at the beginning of this subsection.
\item{(2)}
Suppose that $D$ is a bounded connected $C^2$ open set in $\mathbb{R}^d$ and that $\mu$
denotes the Lebesgue measure on $D$. The reflecting Brownian motion in $\overline{D}$
satisfies the assumptions of the first paragraph in this example and therefore
all the assumptions at the beginning of this subsection.
\item{(3)} Suppose that $D$ is an open subset of  $\mathbb{R}^d$ with finite Lebesgue measure
and that $\mu$ denotes the Lebesgue measure on $D$. Then the subprocesses in $D$ of any of
the subordinate Brownian motions studied in \cite{KSV1, KSV2} satisfy the assumptions of the first paragraph in this example and therefore
all the assumptions at the beginning of this subsection.
\end{description}
}
\end{example}

\subsection{Branching Markov process}\label{subs:ss}

In this subsection, we spell out our assumptions on the branching Markov process.

The branching Markov process $\{X_t: t\ge 0\}$ on $E$ we are going to work with is
determined by three parameters: a spatial motion $\xi=\{\xi_t, \Pi_x\}$ on $E$ satisfying the
assumptions at the beginning of the previous subsection,
 a branching rate function $\beta(x)$ on $E$ which is a non-negative bounded measurable function
and an offspring distribution $\{p_n(x): n=0, 1,, 2, \dots\}$ satisfying the assumption
\begin{equation}\label{1.16}
  \sup_{x\in E}\sum_{n=0}^\infty n^2p_n(x)<\infty.
\end{equation}
We denote the generating function of the offspring distribution by
$$
  \varphi(x,z)=\sum_{n=0}^\infty p_n(x)z^n,\quad x\in E,\quad |z|\leq 1.
$$

Consider a branching system on $E$ characterized by the following properties:
(i) each individual has a random birth and death time;
(ii) given that an individual is born at $x\in E$, the conditional distribution of its
path is determined by $\Pi_x$;
(iii) given the path  $\xi$ of an individual up to time $t$
and given that the particle is alive at time $t$ , its
probability of dying in the interval $[t, t +dt)$ is $\beta(\xi_t)dt + o(dt)$;
(iv) when an individual dies at $x\in E$,
it splits into $n$ individuals all positioned at $x$, with probability $p_n(x)$;
(v) when an individual reaches $\partial$, it disappears from the system;
(vi) all the individuals, once born, evolve independently.

Let $\mathcal{M}_a(E)$ be the space of finite atomic measures on $E$,
and let $\mathcal{B}_b(E)$  be the set of bounded Borel measurable functions on $E$.
Let $X_t(B)$ be the number of particles alive at time $t$ located in $B\in \mathcal{B}(E)$.
Then $X=\{X_t,t\geq 0\}$ is an $\mathcal{M}_a(E)$-valued Markov process.
For any $\nu\in\mathcal{M}_a(E)$, we denote the law of $X$ with initial configuration $\nu$ by $\P_\nu$.
As usual, $\langle f,\nu\rangle:=\int_E f(x)\,\nu(dx)$.
For $0\le f\in \mathcal{B}_b(E)$, let
$$
  \omega(t,x):=\P_{\delta_x}e^{-\langle f,X_t\rangle},
$$
then $\omega(t,x)$ is the unique positive solution to the equation
\begin{equation}\label{1.3}
  \omega(t,x)=\Pi_x\int_0^t \psi(\xi_s,\omega(t-s,\xi_s))\,ds+\Pi_x(e^{-f(\xi_t)}),
\end{equation}
where $\psi(x,z)=\beta(x)(\varphi(x,z)-z),x\in E, z\in [0,1],$ while $\psi(\partial,z)=0, z\in [0,1]$.
By the branching property, we have
$$
\P_{\nu}e^{-\langle f,X_t\rangle}=e^{\langle\log\omega(t,\cdot),\nu\rangle}.
$$
Define
\begin{equation}\label{e:alpha}
\alpha(x):=\beta(x)\left(\sum_{n=1}^\infty np_n(x)-1\right)\quad \mbox{and }
A(x):=\beta(x)\sum_{n=2}^\infty (n-1)n p_n(x).
\end{equation}
By \eqref{1.16}, there exists $K>0$, such that
\begin{equation}\label{1.5}
  \sup_{x\in E}\left(|\alpha(x)|+A(x)\right)\le K.
\end{equation}
For any $f\in\mathcal{B}_b(E)$ and $(t, x)\in (0, \infty)\times E$, define
\begin{equation}\label{1.26}
   T_tf(x):=\Pi_x \left[e^{\int_0^t\alpha(\xi_s)\,ds}f(\xi_t)\right].
\end{equation}
It is well-known that $T_tf(x)=\P_{\delta_x}\langle f,X_t\rangle$ for every $x\in E$.

For any $(t, x,y)\in (0, \infty)\times E\times E$, define
\begin{eqnarray*}
I_0(t, x, y)&=&p_t(x, y),\\
I_n(t, x, y)&=&\int^t_0\int_Ep_s(x, z)I_{n-1}(t-s, z, y)\alpha(z)\mu(dz)ds, \qquad n\ge 1.
\end{eqnarray*}
By induction we can see  that for any $f\in\mathcal{B}_b(E)$ and  any $n\ge 0$,
\begin{equation}\label{e:taylor}
 \int_EI_n(t, x, y)f(y)\mu(dy)
 =\frac1{n!}\Pi_x\left[\left(\int_0^t\alpha(\xi_s)\,ds\right)^nf(\xi_t)\right],
\qquad (t, x)\in (0, \infty)\times E,
\end{equation}
\begin{equation}\label{e:approx}
|I_n(t, x, y)|\le \frac{(\|\alpha\|_\infty t)^n}{n!}p_t(x, y), \qquad (t, x, y)\in (0, \infty)\times E\times E.
\end{equation}
Thus
\begin{equation}\label{e:denrel}
q_t(x, y):=\sum^\infty_{n=0}I_n(t, x, y)\le e^{\|\alpha\|_\infty t}p_t( x, y), \qquad (t, x, y)\in (0, \infty)\times E\times E,
\end{equation}
and, for each $t>0$, the series above converges locally uniformly. Similarly we also have
 $q_t(x, y)\ge \exp(-\|\alpha\|_\infty t)p_t(x, y)$
 for all $(t, x, y)\in (0, \infty)\times E\times E$.
For any $\epsilon\in (0, t/2)$ and $(x, y)\in E\times E$ , we have
\begin{eqnarray*}
\left|\int^\epsilon_0\int_Ep_s(x, z)p_{t-s}(z, y)\alpha(z)\mu(dz)ds\right|
&\le&\|\alpha\|_\infty\int^\epsilon_0\int_Ep_s(x, z)p_{t-s}(z, y)\mu(dz)ds\\
&=&\epsilon\|\alpha\|_\infty p_t(x, y).
\end{eqnarray*}
Similarly, for any $(x, y)\in E\times E$ , we have
$$
\left|\int^t_{t-\epsilon}\int_Ep_s(x, z)p_{t-s}(z, y)\alpha(z)\mu(dz)ds\right|\le \epsilon\|\alpha\|_\infty p_t( x, y).
$$
Hence, for any $t>0$, as $\epsilon\to 0$,
$$
\int^{t-\epsilon}_\epsilon \int_Ep_s(x, z)p_{t-s}(z, y)\alpha(z)\mu(dz)ds
$$
converges to $I_1(t, x, y)$ locally uniformly. For $s\in (\epsilon, t-\epsilon)$ we have,
$$
p_s(x, z)p_{t-s}( z, y)
\le (\widetilde{a}_\epsilon(x)\widetilde{a}_\epsilon(y))^{1/2}
\widetilde{a}_\epsilon(z), \qquad (x, y, z)\in E\times E\times E,
$$
thus it follows from (a) and the dominated convergence theorem that $I_1(t, x, y)$ is continuous on $E\times E$.
Using \eqref{e:approx} and induction, we can show that, for each $n>1$ and $t>0$, $I_n(t, x, y)$ is
continuous on $E\times E$.
$I_1(t, x, y)$ is obviously symmetric in $x$ and $y$.
Using \eqref{e:taylor} and some standard arguments (see the proof of \cite[Theorem 3.10]{CZ}),
one can easily show that, for each $n>1$ and $t>0$, $I_n(t, x, y)$ is symmetric on $E\times E$.
Thus, for any $t>0$, $q_t(x,y)$ is a continuous strictly positive symmetric function on $E\times E$
and for any bounded Borel function $f$ and any $(t, x)\in (0, \infty)\times E$,
$$
  T_tf(x)=\int_E q_t(x,y)f(y)\mu(dy).
$$
It follows immediately from \eqref{e:denrel} that, for any $p\ge 1$, $\{T_t: t\ge 0\}$ is a
strongly continuous semigroup on $L^p(E, \mu)$ and
$$
\|T_tf\|_p^p\le e^{pKt}\|f\|_p^p.
$$

Define $a_t(x):=q_t(x, x)$. It follows from \eqref{e:denrel} and the assumptions (a) and (b) in the
previous subsection that $a_t$ enjoys the following properties.
\begin{description}
  \item[(i)] For any $t>0$, we have
$$
        \int_E a_t(x)\mu(dx)<\infty.
$$
  \item[(ii)] There exists $t_0>0$ such that for all $t\ge t_0$, $a_{t}(x)\in L^2(E,\mu)$.
\end{description}
It follows from (i) above that, for any $t>0$, $T_t$ is a Hilbert-Schmidt operator and thus a compact
operator. Let $L$ be the infinitesimal generator of $\{T_t:t\geq 0\}$ in $L^2(E, \mu)$.
$L$ has purely discrete spectrum with eigenvalues
$-\lambda_1>-\lambda_2>-\lambda_3>\cdots$, and the first eigenvalue $-\lambda_1$ is simple and the eigenfunction
$\phi_1$ associated with $-\lambda_1$ can be chosen to be strictly positive everywhere and continuous.
We will assume that $\|\phi_1\|_2=1$. $\phi_1$ is sometimes denoted as $\phi^{(1)}_1$.
For $k>1$, let $\{\phi^{(k)}_j,j=1,2,\cdots n_k\}$ be
an orthonormal basis of the eigenspace (which is finite dimensional) associated with $-\lambda_k$.
It is well-known that $\{\phi^{(k)}_j,j=1,2,\cdots n_k; k=1,2,\dots\}$ forms a complete orthonormal basis of $L^2(E,\mu)$ and all the eigenfunctions are continuous.
For any $k\ge 1$, $j=1, \dots, n_k$ and $t>0$, we have $T_t\phi^{(k)}_j(x)=e^{-\lambda_k t}\phi^{(k)}_j(x)$ and
\begin{equation}
\label{1.37}
e^{-\lambda_kt/2}|\phi^{(k)}_j|(x)\le a_t(x)^{1/2}, \qquad x\in E.
\end{equation}
It follows from the relation above that all the eigenfunctions $\phi^{(k)}_j$ belong to $L^4(E, \mu)$.
For any $x, y\in E$ and $t>0$, we have
\begin{equation}\label{eq:p}
  q_t(x, y)=\sum^\infty_{k=1}e^{-\lambda_k t}\sum^{n_k}_{j=1}\phi^{(k)}_j(x)\phi^{(k)}_j(y),
\end{equation}
where the series is locally uniformly convergent on $E\times E$.
For the basic facts in the paragraph, one can
refer to  \cite[Section 2]{DS}.

In this paper, we always assume that the branching Markov process $X$ is supercritical,
that is, $\lambda_1<0$.

We will use $\{{\cal F}_t: t\ge0\}$ to denote the filtration of $X$, that is
${\cal F}_t=\sigma(X_s: s\in [0, t])$.
Using the expectation formula of $\langle \phi_1, X_t\rangle$ and the Markov property of $X$,
it is not hard to prove that (see Lemma \ref{lem:1.2} for a proof),
for any nonzero $\nu\in {\cal M}_a(E)$, under $\P_{\nu}$,
the process $W_t:=e^{\lambda_1 t}\langle \phi_1, X_t\rangle$ is a positive martingale. Therefore it converges:
$$
  W_t \to W_\infty,\quad \P_{\nu}\mbox{-a.s.} \quad \mbox{ as }t\to \infty.
$$
Using the assumption \eqref{1.16}
we can show that, as $t\to \infty$, $W_t$ also converges in $L^2(\P_{\nu})$,
so $W_\infty$ is non-degenerate and
 the second moment is finite. Moreover, we have $\P_{\nu}(W_\infty)=\langle\phi_1, \nu\rangle$.
Put $\mathcal{E}=\{W_\infty=0\}$, then $\P_{\nu}(\mathcal{E})<1$.
It is clear that  $\mathcal{E}^c\subset\{X_t(E)>0,\forall t\ge 0\}$.

We will use $\langle\cdot, \cdot\rangle$ to denote inner product in $L^2(E, \mu)$.
Any $f\in L^2(E,\mu)$ admits the following expansion:
\begin{equation}\label{exp:f}
  f(x)=\sum_{k=1}^\infty\sum^{n_k}_{j=1} a_j^k\phi_j^{(k)}(x),
\end{equation}
where $a_j^k=\langle f,\phi_j^{(k)}\rangle$ and the series converges in $L^2(E,\mu)$.
$a^1_1$ will sometimes be written as $a_1$.

\subsection{Main results}

For $f\in L^2(E,\mu)$, define
$$
  \gamma(f):=\inf\{k\geq 1: \mbox{ there exists } j \mbox{ with }
  1\leq j\leq n_k\mbox{ such that }
  a_j^k\neq 0\},
$$
where we use the usual convention $\inf\varnothing=\infty$.
We note that if $f\in L^2(E,\mu)$ is nonnegative and $\mu(x: f(x)>0)>0$,
then $\langle f,\phi_1\rangle>0$ which implies $\gamma(f)=1$.
Define
\begin{eqnarray*}
  f_{(s)}(x)& :=& \sum_{2\lambda_k<\lambda_1}
  \sum_{j=1}^{n_k}a_j^k\phi_j^{(k)}(x),\\
  f_{(c)}(x)&:=&\sum_{2\lambda_k=\lambda_1}
  \sum_{j=1}^{n_k}a_j^k\phi_j^{(k)}(x),\\
  f_{(l)}(x)&:=&f(x)-  f_{(s)}(x)-f_{(c)}(x),\\
  f_1(x)&:=&\sum_{j=1}^{n_{\gamma(f)}}a_j^{\gamma(f)}\phi_j^{(\gamma(f))}(x),\\
  \tilde{f}(x)&:=&f(x)-f_1(x).
\end{eqnarray*}

The main results of this paper are stated in three separate cases:
$\lambda_1>2\lambda_{\gamma(f)}$; $\lambda_1=2\lambda_{\gamma(f)}$
and $\lambda_1<2\lambda_{\gamma(f)}$. When
the branching rate $\beta(x)$ and the offspring distribution $\{p_n(x);
n=0, 1, \dots\}$ are both independent of $x$, the function $\alpha$ defined
in \eqref{e:alpha} reduces to a constant and the eigenvalues $-\lambda_k$ of
$L$ are related to the eigenvalues $-\widetilde{\lambda}_k$ of the generator of $\xi$
by $-\lambda_k=
-\widetilde{\lambda}_k+\alpha$. Therefore,
$\lambda_1>2\lambda_{\gamma(f)}$; $\lambda_1=2\lambda_{\gamma(f)}$
and $\lambda_1<2\lambda_{\gamma(f)}$ are equivalent to
$\alpha>2\widetilde{\lambda}_{\gamma(f)}-\widetilde{\lambda}_1$;
$\alpha=2\widetilde{\lambda}_{\gamma(f)}-\widetilde{\lambda}_1$
and $\alpha<2\widetilde{\lambda}_{\gamma(f)}-\widetilde{\lambda}_1$
respectively. Because of this, when the branching rate $\beta(x)$ and offspring distribution $\{p_n(x);
n=0, 1, \dots\}$ are both independent of $x$,
the cases $\lambda_1>2\lambda_{\gamma(f)}$;
$\lambda_1=2\lambda_{\gamma(f)}$ and $\lambda_1<2\lambda_{\gamma(f)}$
are called the large branching rate case, the critical branching rate case
and the small branching rate case respectively in \cite{RP} and \cite{RSZ}.
Therefore in this paper, even when the branching rate $\beta(x)$ and offspring distribution $\{p_n(x);
n=0, 1, \dots\}$ depend on $x$, we still call the cases
$\lambda_1>2\lambda_{\gamma(f)}$;
$\lambda_1=2\lambda_{\gamma(f)}$ and $\lambda_1<2\lambda_{\gamma(f)}$
the large branching rate case, the critical branching rate case
and the small branching rate case respectively.
Here are the main results of this paper.

\subsubsection{The large branching rate case: $\lambda_1>2\lambda_{\gamma(f)}$}

Define
\begin{equation}\label{1.11}
  H_t^{k,j}:=e^{\lambda_{k}t}\langle\phi_j^{(k)},X_t\rangle.
\end{equation}
$H_t^{1,1}$ will sometimes be written as $H^1_t$.
One can show (see Lemma \ref{lem:1.2} below) that, if $\lambda_1>2\lambda_k$,
then, for any nonzero $\nu\in {\cal M}_a(E)$,
$H_t^{k,j}$ is a martingale under
$\P_{\nu}$
 and bounded in
 $L^2(\P_{\nu})$,
and thus the limit
$ H_\infty^{k,j}:=\lim_{t\to \infty}H_t^{k,j}$ exists
$\P_{\nu}$-a.s. and in $L^2(\P_{\nu})$.

\begin{thrm}\label{The:1.2}
If $f\in L^2(E,\mu)\cap L^4(E,\mu)$ with $\lambda_1>2\lambda_{\gamma(f)}$,
then for any nonzero $\nu\in {\cal M}_a(E)$, as $t\to\infty$,
$$
  e^{\lambda_{\gamma(f)}t}\langle f, X_t\rangle\to \sum_{j=1}^{n_{\gamma(f)}}a^{\gamma(f)}_jH^{\gamma(f),j}_\infty,
   \quad \mbox{ in } L^2(\P_{\nu}).
$$
\end{thrm}

\begin{remark}\label{rem:large}
Suppose $f\in L^2(E,\mu)\cap L^4(E,\mu)$.
When $\gamma(f)=1$, $H_t^1$ reduces to $W_t$, and thus $H_\infty^1=W_\infty$.
Therefore by Theorem~\ref{The:1.2} and the fact that $a_1=\langle f,\phi_1\rangle$,
we get that for any nonzero $\nu\in{\cal M}_a(E)$,
\begin{equation*}
  e^{\lambda_1 t}\langle f,X_t\rangle\to \langle f,\phi_1\rangle W_\infty, \quad \mbox{in }
  L^2(\P_{\nu}),
\end{equation*}
as $t\to\infty$.
In particular, the convergence also holds in $\P_{\nu}$-probability.
\end{remark}

\subsubsection{The small branching rate case: $\lambda_1<2\lambda_{\gamma(f)}$}

Define
\begin{equation}\label{e:sigma}
  \sigma_f^2:=\int_0^\infty e^{\lambda_1 s}\langle A(T_s f)^2,\phi_1\rangle \,ds+\langle f^2,\phi_1\rangle.
\end{equation}

\begin{thrm}\label{The:1.3}
If $f\in L^2(E,\mu)\cap L^4(E,\mu)$ with $\lambda_1<2\lambda_{\gamma(f)}$,
then $ \sigma_f^2<\infty$ and,
for any nonzero $\nu\in {\cal M}_a(E)$, it holds under $\P_{\nu}(\cdot\mid \mathcal{E}^c)$ that
$$
  \left(e^{\lambda_1 t}\langle \phi_1, X_t\rangle, ~\frac{\langle f , X_t\rangle}{\sqrt{\langle \phi_1,X_t\rangle}} \right)\stackrel{d}{\rightarrow}(W^*,~G_1(f)), \quad t\to\infty,
$$
where $W^*$ has the same distribution as $W_\infty$ conditioned on $\mathcal{E}^c$
and $G_1(f)\sim \mathcal{N}(0,\sigma_f^2)$. Moreover, $W^*$ and $G_1(f)$ are independent.
\end{thrm}

\subsubsection{The critical branching rate case: $\lambda_1=2\lambda_{\gamma(f)}$}

Define
\begin{equation}\label{e:rho}
\rho_f^2:=\left\langle Af_1^2,\phi_1\right\rangle.
\end{equation}

\begin{thrm}\label{The:1.4}
If $f\in L^2(E,\mu)\cap L^4(E,\mu)$ with $\lambda_1=2\lambda_{\gamma(f)}$,
then $ \rho_f^2<\infty$ and,
for any nonzero $\nu \in {\cal M}_a(E)$, it holds under $\P_{\nu}(\cdot\mid \mathcal{E}^c)$ that
$$
  \left(e^{\lambda_1 t}\langle \phi_1, X_t\rangle, ~\frac{\langle f , X_t\rangle}{\sqrt{t\langle \phi_1,X_t\rangle}} \right)\stackrel{d}{\rightarrow}(W^*,~G_2(f)), \quad t\to\infty,
$$
where $W^*$ has the same distribution as $W_\infty$ conditioned on $\mathcal{E}^c$
and $G_2(f)\sim \mathcal{N}(0,\rho_f^2)$. Moreover, $W^*$ and $G_2(f)$ are independent.
\end{thrm}

\subsubsection{Further results in the large branching rate case}

In this subsection we give two central limit theorems for the case
$\lambda_1>2\lambda_{\gamma(f)}$.
Define
\begin{equation}\label{1.27}
   H_\infty:=\sum_{2\lambda_k<\lambda_1}\sum_{j=1}^{n_k}a_j^kH^{k,j}_\infty.
\end{equation}
Let
\begin{equation}\label{1.61}
  \beta_{f}^2:=\int_0^\infty e^{-\lambda_1 s}
    \left\langle A(\sum_{2\lambda_k<\lambda_1}\sum_{j=1}^{n_k}e^{\lambda_ks}a_j^k\phi_j^k)^2,\phi_1\right\rangle\,ds
    -\langle (f_{(s)})^2,\phi_1\rangle.
\end{equation}
In Section 3.3 we will see that
$\beta_{f}^2=\langle Var_{\delta_{\cdot}}H_\infty,~\phi_1\rangle$.

\begin{thrm}\label{The:2.1}
If $f\in L^2(E,\mu)\cap L^4(E,\mu)$ satisfies $\lambda_1>2\lambda_{\gamma(f)}$
and $f_{(c)}=0$, then $\sigma_{f_{(l)}}^2<\infty$ and $\beta_f^2<\infty$.
For any nonzero $\nu\in {\cal M}_a(E)$, it holds under $\P_{\nu}(\cdot\mid \mathcal{E}^c)$ that, as $t\to\infty$,
$$
  \left(e^{\lambda_1 t}\langle \phi_1,X_t\rangle,~\langle \phi_1,X_t\rangle^{-1/2}
  \left(\langle f,X_t\rangle-\sum_{2\lambda_k<\lambda_1}
  e^{-\lambda_kt}\sum_{j=1}^{n_k}a_j^kH^{k,j}_\infty\right) \right)
  \stackrel{d}{\rightarrow}(W^*,~G_3(f)),
$$
where $W^*$ has the same distribution as $W_\infty$ conditioned on $\mathcal{E}^c$,
and $G_3(f)\sim \mathcal{N}(0,\sigma_{f_{(l)}}^2+\beta_{f}^2)$.
Moreover, $W^*$ and $G_3(f)$ are independent.
\end{thrm}

\begin{remark}
If $2\lambda_k<\lambda_1$,
then, for any nonzero $\nu\in{\cal M}_a(E)$, it holds under $\P_{\nu}(\cdot\mid \mathcal{E}^c)$ that,
as $t\to\infty$,
$$
  \left(e^{\lambda_1 t}\langle \phi_1,X_t\rangle,~\frac{\left(\langle \phi_j^{(k)},X_t\rangle-e^{-\lambda_kt}H^{k,j}_\infty\right)}{\langle \phi_1,X_t\rangle^{1/2}} \right)\stackrel{d}{\rightarrow}(W^*,~G_3),
$$
where $G_3\sim\mathcal{N}\left(0,\frac{1}{\lambda_1-2\lambda_k}\langle A(\phi_j^{(k)})^2,\phi_1\rangle\right)$.
In particular, for $\phi_1$, we have
\begin{equation*}
  \left(e^{\lambda_1 t}\langle \phi,X_t\rangle,~\frac{\left(\langle \phi_1,X_t\rangle-e^{-\lambda_1t}W_\infty\right)}{\langle \phi_1,X_t\rangle^{1/2}}\right)\stackrel{d}{\rightarrow}(W^*,~G_3), \quad t\to\infty,
\end{equation*}
where $G_3\sim\mathcal{N}\left(0,-\frac{1}{\lambda_1}\int_E A(x)(\phi_1(x))^3\mu(dx)\right)$.
\end{remark}

\begin{thrm}\label{The:2.3}
If $f\in L^2(E,\mu)\cap L^4(E,\mu)$ satisfies
$\lambda_1>2\lambda_{\gamma(f)}$
 and $f_{(c)}\ne0$,
then, for any nonzero $\nu\in{\cal M}_a(E)$, it holds under $\P_{\nu}(\cdot\mid \mathcal{E}^c)$ that,
as $t\to\infty$,
$$
  \left(e^{\lambda_1 t}\langle \phi_1,X_t\rangle,~t^{-1/2}\langle \phi_1,X_t\rangle^{-1/2}
  \left(\langle f,X_t\rangle-\sum_{\lambda_k<\lambda_1/2}e^{-\lambda_kt}\sum_{j=1}^{n_k}a_j^kH^{k,j}_\infty\right)\right)
  \stackrel{d}{\rightarrow}(W^*,~G_4(f)),
$$
where $W^*$ has the same distribution as $W_\infty$ conditioned on $\mathcal{E}^c$,
and $G_4(f)\sim \mathcal{N}(0,\rho_{f_{(c)}}^2)$.
Moreover, $W^*$ and $G_4(f)$ are independent.
\end{thrm}

\begin{remark} {\rm By combining the techniques of this paper with the backbone decomposition
of superprocesses (see \cite{BKS}),
one can extend the central limit theorems, for super-OU processes, of \cite{RSZ} to
superprocesses with spatial-dependent branching mechanisms and with
spatial motions satisfying the assumptions (a) and (b).}
\end{remark}

\section{Preliminaries}

In this section, we will give the estimates on the moments of the branching Markov process $X$.

\subsection{Estimates on the semigroup $T_t$}

In the remainder of this paper we will use the  following notation:
for two positive functions $f$ and $g$ on $E$, $f(x)\lesssim g(x)$ for $x\in E$ means that there exists a constant $c>0$ such that
$f(x)\le cg(x)$ for any $x\in E$.

\begin{lemma}\label{lem:expansion}
 For any $f\in L^2(E,\mu)$, $x\in E$ and $t>0$, we have
\begin{equation}\label{1.17}
T_t f(x)=\sum_{k=\gamma(f)}^\infty e^{-\lambda_k t}\sum_{j=1}^{n_k}a_j^k\phi^{(k)}_j(x)
\end{equation}
and
\begin{equation}\label{1.25}
  \lim_{t\to\infty}e^{\lambda_{\gamma(f)}t}T_t f(x)=\sum_{j=1}^{n_{\gamma(f)}}a_j^{\gamma(f)}\phi_j^{(\gamma(f))}(x),
\end{equation}
where the series in \eqref{1.17} converges absolutely and uniformly in any compact subset of $E$.
Moreover, for any $t_1>0$,
\begin{eqnarray}
   &&\sup_{t>t_1}e^{\lambda_{\gamma(f)}t}|T_t f(x)|\lesssim (a_{t_1}(x))^{1/2},\label{1.36}\\
  &&\sup_{t>t_1}e^{(\lambda_{\gamma(f)+1}-\lambda_{\gamma(f)})t}
  \left|e^{\lambda_{\gamma(f)}t}T_t f(x)-f_1(x)\right|
  \lesssim (a_{t_1}(x))^{1/2}.\label{1.43}
\end{eqnarray}
\end{lemma}

\textbf{Proof:}
 Using \eqref{eq:p}, it is easy to see that for any
 $(t, x)\in (0, \infty)\times E$,
 $$
 T_t f(x)=\int_E \sum^\infty_{k=1}e^{-\lambda_k t}\sum^{n_k}_{j=1}\phi^{(k)}_j(x)\phi^{(k)}_j(y)f(y)\,\mu(dy).
 $$
To prove \eqref{1.17}, we only need to show that, for any $t_1>0$ and any $(t, x)\in (t_1, \infty)\times E$,
 $$
 \sum^\infty_{k=1}e^{-\lambda_k t}\sum^{n_k}_{j=1}|\phi^{(k)}_j(x)|\int_E|\phi^{(k)}_j(y)||f(y)|\,\mu(dy)<\infty,
 $$
 and that the series convergent uniformly on any compact subset of $E$.
 By H\"{o}lder's inequality, we get $\int_E|\phi^{(k)}_j(y)||f(y)|\,\mu(dy)\le \|f\|_2$.
 Then by \eqref{1.37}, for $(t, x)\in (t_1, \infty)\times E$, we have
\begin{eqnarray}\label{1.39}
  &&\sum^\infty_{k=1}e^{-\lambda_k t}\sum^{n_k}_{j=1}|\phi^{(k)}_j(x)|\int_E|\phi^{(k)}_j(y)||f(y)|\,\mu(dy)
     \le  \sum^\infty_{k=1}n_ke^{-\lambda_k (t-t_1/2)}\|f\|_2a_{t_1}(x)^{1/2} \\
   &\le &e^{-\lambda_1(t-t_1)} \|f\|_2a_{t_1}(x)^{1/2} \sum^\infty_{k=1}n_ke^{-\lambda_k t_1/2}.
\end{eqnarray}
By \eqref{eq:p}, we have
\begin{equation}\label{1.40}
   \sum_{k=1}^\infty e^{-\lambda_k t_1/2}\sum_{j=1}^{n_k}|\phi^{(k)}_j(x)|^2=a_{t_1/2}(x), \qquad x\in E.
\end{equation}
Consequently, integrating both sides of \eqref{1.40},
\begin{equation}\label{1.41}
  \sum_{k=1}^\infty n_k e^{-\lambda_k t_1/2}= \int_E a_{t_1/2}(x)\mu(dx)<\infty.
\end{equation}
Thus, for any $(t, x)\in (t_1, \infty)\times E$,
\begin{equation}\label{1.421}
  \sum^\infty_{k=1}e^{-\lambda_k t}\sum^{n_k}_{j=1}|\phi^{(k)}_j(x)|\int_E|\phi^{(k)}_j(y)||f(y)|\,\mu(dy)
  \le e^{-\lambda_1(t-t_1)}\|f\|_2a_{t_1}(x)^{1/2}\int_E a_{t_1/2}(x)\mu(dx).
\end{equation}
 Thus, we get \eqref{1.17}.
 The above argument shows that for any $(t, x)\in (t_1, \infty)\times E$,
 \begin{equation}\label{1.42}
   |T_t f(x)|\le \sum^\infty_{k=\gamma(f)}e^{-\lambda_k t}\sum^{n_k}_{j=1}|a_j^k||\phi^{(k)}_j(x)|
   \le e^{-\lambda_{\gamma(f)}(t-t_1)}\|f\|_2a_{t_1}(x)^{1/2}\int_E a_{t_1/2}(x)\mu(dx),
 \end{equation}
 which implies \eqref{1.36}.

Applying \eqref{1.36} to $\tilde{f}$, we obtain that for $x\in E$,
$$
\sup_{t>t_1}e^{\lambda_{\gamma(f)+1}t}|T_t(\tilde{f})(x)|\lesssim (a_{t_1}(x))^{1/2}.
$$
Now  \eqref{1.25} and \eqref{1.43} follow immediately.\hfill$\Box$

The proof of the lemma above also yields the following result which will be used later.
\begin{lemma}\label{lem:rsnew}
Suppose that $\{f_t(x):t>0\}$ is a family of functions in $L^2(E, \mu)$.
 If $\lim_{t\to\infty} \|f_t\|_2=0$, then for any $x\in E$,
$$
\lim_{t\to\infty}e^{\lambda_1t}T_t f_t(x)=0.
$$
\end{lemma}

\textbf{Proof:}
Applying \eqref{1.42} to $f_t$ and using the fact $\lambda_1\le \lambda_{\gamma(f_t)}$, we get that
for any $(t, x)\in (t_0, \infty)\times E$,
\begin{eqnarray*}
|T_t f_t(x)|\le e^{-\lambda_1(t-t_0)}\|f_t\|_2(a_{t_0}(x))^{1/2}\int_E a_{t_0/2}(x)\mu(dx).
\end{eqnarray*}
Thus, for any $(t, x)\in (t_0, \infty)\times E$,
$$
|e^{\lambda_1t}T_t f_t(x)|\le e^{\lambda_1t_0}(a_{t_0}(x))^{1/2}\|f_t\|_2\int_E a_{t_0/2}(x)\mu(dx),
$$
from which the assertion of the lemma follows immediately.
\hfill$\Box$

\subsection{Estimates on the second moment of the branching Markov process}

Recall the formula for the second moment of the branching Markov process $\{X_t: t\ge 0\}$
(see, for example, \cite[Lemma 3.3]{Sh}):
for $f\in \mathcal{B}_b(E)$, we have for any $(t, x)\in (0, \infty)\times E$,
\begin{equation}\label{1.13}
   \P_{\delta_x}\langle f,X_t\rangle^2=\int_0^tT_{s}[A(T_{t-s}f)^2](x)\,ds+ T_t(f^2)(x).
\end{equation}
For any $f\in L^2(E,\mu)\cap L^4(E,\mu)$ and $x\in E$, by $(T_{t-s}f)^2(x)\le e^{K(t-s)}T_{t-s}(f^2)(x)$, we have
\begin{equation*}
  \int_0^tT_{s}[A(T_{t-s}f)^2](x)\,ds\le Ke^{K(t-s)}tT_t(f^2)(x)<\infty,
\end{equation*}
which implies $\int_0^tT_{s}[A(T_{t-s}f)^2](x)\,ds+ T_t(f^2)(x)<\infty$.
Thus, using a routine limit argument, one can easily check
that \eqref{1.13} also holds for $f\in L^2(E,\mu)\cap L^4(E,\mu)$.

\begin{lemma} \label{lem:2.2}
Assume that $f\in L^2(E,\mu)\cap L^4(E,\mu)$.
\begin{description}
  \item{(1)}
   If $\lambda_1<2\lambda_{\gamma(f)}$, then for any $x\in E$,
    \begin{equation}\label{limit-mean}
   \lim_{t\to \infty}e^{\lambda_1t/2}\P_{\delta_x}\langle f, X_t\rangle=0,
   \end{equation}
  \begin{equation}\label{2.2}
   \lim_{t\to \infty}e^{\lambda_1t}\P_{\delta_x}\langle f,X_t\rangle^2
   = \int_0^\infty e^{\lambda_1 s}
   \langle A(T_s f)^2,\phi_1\rangle \,ds \phi_1(x)+\langle f^2,\phi_1\rangle \phi_1(x).
  \end{equation}
  Moreover, for $(t, x)\in (3t_0, \infty)\times E$, we have
  \begin{equation}\label{2.10}
     e^{\lambda_1t}Var_{\delta_x}\langle f,X_t\rangle\lesssim a_{t_0}(x)^{1/2}.
  \end{equation}

  \item{(2)}
  If $\lambda_1=2\lambda_{\gamma(f)}$, then for any $(t, x)\in (3t_0, \infty)\times E$,
  \begin{equation}\label{1.49}
   \left|t^{-1}e^{\lambda_1 t}Var_{\delta_x}\langle f,X_t\rangle-\rho_f^2\phi_1(x)\right|
   \lesssim t^{-1}\left(a_{t_0}(x)^{1/2}+a_{t_0}(x)\right),
  \end{equation}
  where $\rho^2_f$ is defined by \eqref{e:rho}.

  \item{(3)}
If $\lambda_1>2\lambda_{\gamma(f)}$, then for any $x\in E$,
  \begin{equation}\label{1.23}
    \lim_{t\to \infty}e^{2\lambda_{\gamma(f)}t}\P_{\delta_x}\langle f,X_t\rangle^2=\eta_f^2(x),
  \end{equation}
  where
  \begin{equation*}\label{e:etaf}
   \eta_f^2(x):=\int_0^\infty e^{2\lambda_{\gamma(f)}s}T_s(Af^2_1)(x)\,ds.
  \end{equation*}
  Moreover,
  for any $(t, x)\in (3t_0, \infty)\times E$,
  \begin{equation}\label{1.60}
   e^{2\lambda_{\gamma(f)}t}\P_{\delta_x}\langle f,X_t\rangle^2\lesssim a_{t_0}(x)^{1/2}.
  \end{equation}
\end{description}
\end{lemma}

\textbf{Proof:}\quad
(1)
If $\lambda_1<2\lambda_{\gamma(f)}$, then
by \eqref{1.26} and \eqref{1.36}, we have for any $(t, x)\in (t_0, \infty)\times E$,
\begin{eqnarray}\label{1.50}
e^{\lambda_1t/2}\left|\P_{\delta_x}\langle f,X_t\rangle\right|&=&e^{(\lambda_1-2\lambda_{\gamma(f)})t/2}[e^{\lambda_{\gamma(f)}t}|T_tf(x)|]\nonumber\\
&\lesssim& e^{(\lambda_1-2\lambda_{\gamma(f)})t/2} a_{t_0}(x)^{1/2}\to 0,\quad\mbox{as }t\to\infty.
\end{eqnarray}
In the remainder of the proof of (1), we always assume $t>3t_0$. It follows from \eqref{1.13} that for any $x\in E$,
\begin{eqnarray}\label{1.20}
   &&e^{\lambda_1t}\P_{\delta_x}\langle f,X_t\rangle^2=
   e^{\lambda_1t}\int^t_0T_{t-s}[A(T_sf)^2](x)\,ds+ e^{\lambda_1t}T_t(f^2)(x)\nonumber\\&=&
   \int_0^t e^{(\lambda_1-2\lambda_{\gamma(f)})s}e^{\lambda_1(t-s)}T_{t-s}
       [A(e^{\lambda_{\gamma(f)}s}T_s f)^2](x)\,ds+ e^{\lambda_1t}T_t(f^2)(x)\nonumber\\
   &=&\left(\int_0^{t-t_0}+\int_{t-t_0}^t\right)e^{(\lambda_1-2\lambda_{\gamma(f)})s}e^{\lambda_1(t-s)}T_{t-s}
       [A(e^{\lambda_{\gamma(f)}s}T_s f)^2](x)\,ds+ e^{\lambda_1t}T_t(f^2)(x)\nonumber\\
   &=:&V_1(t,x)+V_2(t,x)+e^{\lambda_1t}T_t(f^2)(x).
\end{eqnarray}
For $V_1(t,x)$, we claim that for $s<t-t_0$, we have
\begin{equation}\label{1.14}
   e^{\lambda_1(t-s)}T_{t-s}[A(e^{\lambda_{\gamma(f)}s}T_s f)^2](x)\lesssim a_{t_0}(x)^{1/2}, \qquad x\in E.
\end{equation}
If $s\leq t_0$, using $(T_sf(x))^2\leq e^{Ks}T_s(f^2)(x)$ and
\eqref{1.36}, we obtain that for any $x\in E$,
$$
e^{\lambda_1(t-s)}T_{t-s}[A(e^{\lambda_{\gamma(f)}s}T_s f)^2](x)
\leq Ke^{Ks}e^{-(\lambda_1-2\lambda_{\gamma(f)})s}e^{\lambda_1t}T_t(f^2)(x)\lesssim a_{t_0}(x)^{1/2}.
$$
If $t_0<s<t-t_0$, by \eqref{1.36}, we have for any $x\in E$,
  \begin{equation}\label{1.29}
   e^{\lambda_1(t-s)}T_{t-s}[A(e^{\lambda_{\gamma(f)}s}T_s f)^2](x)
   \lesssim e^{\lambda_1(t-s)}T_{t-s}a_{t_0}(x)\lesssim  a_{t_0}(x)^{1/2}.
\end{equation}
Thus, we have proved the claim.
By \eqref{1.14}, we get that for any $x\in E$,
\begin{equation}\label{1.15}
  V_1(t,x)\lesssim (a_{t_0}(x))^{1/2}.
\end{equation}
By \eqref{1.25} and the dominated convergence theorem, we easily get that for any $x\in E$,
\begin{equation}\label{1.45}
  \lim_{t\to\infty}V_1(t,x)=\int_0^{\infty}
  e^{\lambda_1s}\langle A(T_s f)^2,\phi_1\rangle \,ds \phi_1(x).
\end{equation}

Now we deal with $V_2(t,x)$. It follows from \eqref{1.36} that $e^{\lambda_{\gamma(f)}s}T_s f(x)\lesssim a_{2t_0}(x)^{1/2}$ for $(s,x)\in (2t_0,\infty)\times E$.
 Thus,
\begin{eqnarray}\label{2.46}
 V_2(t,x)&\lesssim& \int_{t-t_0}^t e^{(\lambda_1-2\lambda_{\gamma(f)})s}e^{\lambda_1(t-s)}T_{t-s}(a_{2t_0})(x)\,ds\nonumber\\
  &=&e^{(\lambda_1-2\lambda_{\gamma(f)})t}\int_{0}^{t_0} e^{2\lambda_{\gamma(f)}s}T_s(a_{2t_0})(x)\,ds \nonumber\\
  &\lesssim& e^{(\lambda_1-2\lambda_{\gamma(f)})t}\int_{0}^{t_0}T_s(a_{2t_0})(x)\,ds.
\end{eqnarray}
  We now show that for any $x\in E$, $\int_{0}^{t_0}T_s(a_{2t_0})(x)\,ds<\infty$. By \eqref{eq:p}, we have
\begin{eqnarray*}
  a_{2t_0}(x)&=&\sum_{k=1}^\infty \sum_{j=1}^{n_k}e^{-2\lambda_k t_0}\left|\phi_j^{(k)}(x)\right|^2=\sum_{k=1}^\infty \sum_{j=1}^{n_k}e^{-\lambda_k t_0}\left|T_{t_0/2}\phi_j^{(k)}(x)\right|^2 \\
  &\le& e^{Kt_0/2}\sum_{k=1}^\infty \sum_{j=1}^{n_k}e^{-\lambda_k t_0}T_{t_0/2}\left|\phi_j^{(k)}\right|^2(x)=e^{Kt_0/2}T_{t_0/2}(a_{t_0})(x).
\end{eqnarray*}
So, by H\"{o}lder's inequality, we have
\begin{eqnarray*}
  T_s(a_{2t_0})(x) &\le& e^{Kt_0/2}T_{s+t_0/2}(a_{t_0})(x)\le e^{Kt_0/2}\|a_{t_0}\|_2a_{2s+t_0}(x)^{1/2}.
\end{eqnarray*}
By \eqref{eq:p}, we have
$$
a_{2s+t_0}(x)=\sum_{k=1}^\infty \sum_{j=1}^{n_k}e^{-\lambda_k (t_0+2s)}\left|\phi_j^{(k)}(x)\right|^2\le e^{-2\lambda_1s}a_{t_0}(x),
$$
which implies
\begin{equation}\label{1.53}
  \int_{0}^{t_0}T_s(a_{2t_0})(x)\,ds\lesssim a_{t_0}(x)^{1/2}.
\end{equation}
Hence for any $x\in E$, as $t\to\infty$,
\begin{equation}\label{1.54}
  V_2(t,x)\lesssim e^{(\lambda_1-2\lambda_{\gamma(f)})t}a_{t_0}(x)^{1/2}\to 0.
\end{equation}
 Thus, by \eqref{1.45} and \eqref{1.54}, we get that for any $x\in E$,
 \begin{equation}\label{2.47}
   \lim_{t\to\infty}e^{\lambda_1t}\int_{0}^{t} T_{t-s}[A(T_s f)^2](x)\,ds
    =\int_{0}^\infty e^{\lambda_1 s}\langle A(T_s f)^2,\phi_1\rangle \,ds \phi_1(x).
 \end{equation}
Since $f^2\in L^2(E,\mu)$, by \eqref{1.25}, we easily get
$\lim_{t\to\infty}e^{\lambda_1t}T_t(f^2)(x)=\langle f^2,\phi_1\rangle\phi_1(x)$ for every $x\in E$,
which implies \eqref{2.2}.

By \eqref{1.36}, we also have $e^{\lambda_1t}T_t(f^2)(x)\lesssim a_{t_0}(x)^{1/2}$ for any $x\in E$.
Combining \eqref{1.15} and \eqref{1.54}, we get that for any $(t,x)\in (3t_0,\infty)\times E$,
$$
e^{\lambda_1t}Var_{\delta_x}\langle f,X_t\rangle\leq e^{\lambda_1t}\P_{\delta_x}\langle f,X_t\rangle^2
\lesssim a_{t_0}(x)^{1/2}.
$$
The proof of (1) is now complete.

(2)
  If $2\lambda_{\gamma(f)}=\lambda_1$,
then by \eqref{1.26} and \eqref{1.13}, we have for any $(t, x)\in (0, \infty)\times E$,
\begin{eqnarray}\label{1.21}
t^{-1}e^{\lambda_1 t}Var_{\delta_x}\langle f,X_t\rangle
&=&t^{-1}\int_0^t e^{\lambda_1s}T_s[A(e^{\lambda_{\gamma(f)}(t-s)}T_{t-s}f)^2](x)\,ds\nonumber\\
&&+ t^{-1}e^{\lambda_1t}T_t(f^2)(x)-t^{-1}\left(e^{\lambda_{\gamma(f)}t}T_tf(x)\right)^2.
\end{eqnarray}
Thus,
\begin{eqnarray}\label{1.6}
  &&\left|t^{-1}e^{\lambda_1 t}Var_{\delta_x}\langle f,X_t\rangle-\rho_f^2\phi_1(x)\right|\nonumber\\
  &\leq& t^{-1}\int_0^{t} e^{\lambda_1s}T_s\left[A\left|(e^{\lambda_{\gamma(f)}(t-s)}T_{t-s}f)^2-f_1^2\right|\right](x)\,ds\nonumber\\
  &&+t^{-1}\int_0^{t} \left|e^{\lambda_1s}T_s(Af_1^2)(x)-\langle Af_1^2,\phi_1\rangle \phi_1(x)\right|\,ds\nonumber\\
   &&+t^{-1}e^{\lambda_1t}T_t(f^2)(x)+t^{-1}\left(e^{\lambda_{\gamma(f)}t}T_tf(x)\right)^2\nonumber\\
   &=:& A_1(t,x)+A_2(t,x)+A_3(t,x)+A_4(t,x).
\end{eqnarray}
In the remainder of the proof of (2), we always assume $t>3t_0$.

For $A_1(t,x)$, by \eqref{1.37}, \eqref{1.36} and \eqref{1.43},
for $t-s>2t_0$,
we have for any $x\in E$,
\begin{eqnarray*}
  \left|(e^{\lambda_{\gamma(f)}(t-s)}T_{t-s}f)^2-f_1(x)^2\right|&\leq&
   \left|e^{\lambda_{\gamma(f)}(t-s)}T_{t-s}f(x)-f_1(x)\right|
  (e^{\lambda_{\gamma(f)}(t-s)}|T_{t-s}f(x)|+|f_1(x)|)\nonumber\\
  &\lesssim & e^{(\lambda_{\gamma(f)}-\lambda_{\gamma(f)+1})(t-s)}a_{2t_0}(x).
  \end{eqnarray*}
So by \eqref{1.36} and \eqref{1.53}, we have for any $(t, x)\in (3t_0, \infty)\times E$,
\begin{eqnarray}\label{1.19}
   &&t^{-1}\int_{0}^{t-2t_0} e^{\lambda_1s}T_s[A|(e^{\lambda_{\gamma(f)}(t-s)}T_{t-s}f)^2-f_1^2|](x)\,ds\\
  &\lesssim& t^{-1}\int_0^{t-2t_0} e^{(\lambda_{\gamma(f)}-\lambda_{\gamma(f)+1})(t-s)} e^{\lambda_1s}T_s(a_{2t_0})(x)\,ds\nonumber\\
  &\lesssim& t^{-1}\int_0^{t_0} e^{(\lambda_{\gamma(f)}-\lambda_{\gamma(f)+1})(t-s)} e^{\lambda_1s}T_s(a_{2t_0})(x)\,ds+ t^{-1}\int_{t_0}^{t-2t_0} e^{(\lambda_{\gamma(f)}-\lambda_{\gamma(f)+1})(t-s)}\,ds a_{t_0}(x)^{1/2}\nonumber\\
  &\lesssim& t^{-1}\int_0^{t_0} T_s(a_{2t_0})(x)\,ds+t^{-1}a_{t_0}(x)^{1/2}\lesssim t^{-1}a_{t_0}(x)^{1/2}.
\end{eqnarray}
Using $(T_{t-s}f(x))^2\leq e^{K(t-s)}T_{t-s}(f^2)(x)$ and \eqref{1.36}, we get that for any $x\in E$,
\begin{eqnarray}
   &&t^{-1}\int_{t-2t_0}^te^{\lambda_1s}T_s\left[A|(e^{\lambda_{\gamma(f)}(t-s)}T_{t-s}f)^2-f_1^2|\right](x)\,ds\nonumber\\
  &\leq& Kt^{-1}\int_{t-2t_0}^te^{\lambda_1s}T_s\left(e^{(2\lambda_{\gamma(f)}+K)(t-s)}
  T_{t-s}(f^2)+f_1^2\right)(x)\,ds \nonumber\\
   &=& Kt^{-1}\int_{t-2t_0}^te^{\lambda_1s}e^{(2\lambda_{\gamma(f)}+K)(t-s)}\,dsT_t(f^2)(x)+
   Kt^{-1}\int_{t-2t_0}^te^{\lambda_1s}T_s(f_1^2)(x)\,ds\nonumber\\
   &=&Kt^{-1}\int_0^{2t_0}e^{-\lambda_1s}e^{(2\lambda_{\gamma(f)}+K)s}\,ds~e^{\lambda_1t}T_t(f^2)(x)+
   Kt^{-1}\int_{t-2t_0}^te^{\lambda_1s}T_s(f_1^2)(x)\,ds\nonumber\\
   &\lesssim&  t^{-1}a_{t_0}(x)^{1/2}\label{1.8}.
\end{eqnarray}
Thus we get that for any $x\in E$,
\begin{equation}\label{1.46}
  A_1(t,x)\lesssim t^{-1}a_{t_0}(x)^{1/2}.
\end{equation}

Next we consider $A_2(x,t)$.
By \eqref{1.43}, we have for $(s, x)\in (t_0, \infty)\times E$,
$$
\left|e^{\lambda_1s}T_s(Af_1^2)(x)-\langle Af_1^2,\phi_1\rangle \phi_1(x)\right|
\lesssim e^{-(\lambda_2-\lambda_1)s}a_{t_0}(x)^{1/2},
$$
which implies
\begin{eqnarray}
 && t^{-1}\int_{t_0}^{t} \left|e^{\lambda_1s}T_s(Af_1^2)(x)-\langle Af_1^2,\phi_1\rangle \phi_1(x)\right|\,ds\nonumber\\
 & \lesssim& t^{-1}\int_{t_0}^{t} e^{-(\lambda_2-\lambda_1)s}\,ds a_{t_0}(x)^{1/2}\lesssim t^{-1}a_{t_0}(x)^{1/2}.
 \label{rsnew}
\end{eqnarray}
By \eqref{1.37}, we get $\phi_1(x)\lesssim a_{t_0}(x)^{1/2}$ and $|f_1(x)|\lesssim a_{2t_0}(x)^{1/2}$
for any $x\in E$. So for any $x\in E$,
\begin{eqnarray*}
  &&t^{-1}\int_0^{t_0}\left|e^{\lambda_1s}T_s(Af_1^2)(x)-\langle Af_1^2,\phi_1\rangle \phi_1(x)\right|\,ds\\
  &\lesssim& Kt^{-1}\int_0^{t_0}e^{\lambda_1s}T_s(a_{2t_0})(x)\,ds+ Kt^{-1}\langle f_1^2,\phi_1\rangle \phi_1(x)\\
  &\lesssim& t^{-1}\int_0^{t_0}T_s(a_{2t_0})(x)\,ds+ t^{-1}a_{t_0}(x)^{1/2}\lesssim t^{-1}a_{t_0}(x)^{1/2}.
\end{eqnarray*}
Thus, we get that for any $x\in E$,
\begin{equation}\label{2.13}
  A_2(t,x)\lesssim t^{-1}a_{t_0}(x)^{1/2}.
\end{equation}

By \eqref{1.36}, we easily get that for any $x\in E$,
$$
A_3(t,x)\lesssim t^{-1}a_{t_0}(x)^{1/2}\quad \mbox{and}\quad A_4(t,x)\lesssim t^{-1}a_{t_0}(x).
$$
Consequently, we have
\begin{equation}\label{1.48}
  \left|t^{-1}e^{\lambda_1 t}Var_{\delta_x}\langle f,X_t\rangle-\rho_f^2\phi_1(x)\right|
  \lesssim t^{-1}\left(a_{t_0}(x)+a_{t_0}(x)^{1/2}\right).
\end{equation}
The proof of (2) is now complete.

(3)
 If $\lambda_1>2\lambda_{\gamma(f)}$,
then by \eqref{1.13}, we have, for $t>3t_0$ and $x\in E$,
\begin{eqnarray}\label{1.31}
&&e^{2\lambda_{\gamma(f)}t} \P_{\delta_x}\langle f,X_t\rangle^2 \nonumber\\
&=&\int_0^te^{-(\lambda_1-2\lambda_{\gamma(f)})s}e^{\lambda_1s}T_s[A(e^{\lambda_{\gamma(f)}(t-s)}T_{t-s} f)^2](x)\,ds
+e^{-(\lambda_1-2\lambda_{\gamma(f)})t}e^{\lambda_1 t}T_t(f^2)(x)\nonumber\\
&=&\left(\int_0^{t_0}+\int_{t_0}^t\right)e^{-(\lambda_1-2\lambda_{\gamma(f)})s}
e^{\lambda_1s}T_s[A(e^{\lambda_{\gamma(f)}(t-s)}T_{t-s} f)^2](x)\,ds
+e^{-(\lambda_1-2\lambda_{\gamma(f)})t}e^{\lambda_1 t}T_t(f^2)(x)\nonumber\\
&=&B_1(t,x)+B_2(t,x)+B_3(t,x).
\end{eqnarray}
In the remainder of this proof, we always assume $t>3t_0$.
For $s\leq t_0$, we get $t-s>2t_0$. So for any $x\in E$,
$$
(e^{\lambda_{\gamma(f)}(t-s)}T_{t-s} f(x))^2\lesssim a_{2t_0}(x),
$$
which implies
$$
B_1(t,x)\lesssim \int_0^{t_0}e^{-(\lambda_1-2\lambda_{\gamma(f)})s}e^{\lambda_1s}T_s(a_{2t_0})(x)\,ds\lesssim\int_0^{t_0}T_s(a_{2t_0})(x)\,ds
\lesssim a_{t_0}(x)^{1/2}.
$$
Thus by the dominated convergence theorem, we get that for any $x\in E$,
\begin{equation*}
  \lim_{t\to\infty}B_1(t,x)
  =\int_{0}^{t_0} e^{-(\lambda_1-2\lambda_{\gamma(f)})s}e^{\lambda_1s}
 T_s(Af^2_1)(x)\,ds.
\end{equation*}
Now we consider $B_2(t,x)$. Using \eqref{1.14}, we get, for $(s, x)\in (t_0, \infty)\times E$,
$$
e^{\lambda_1s}T_s[A(e^{\lambda_{\gamma(f)}(t-s)}T_{t-s} f)^2](x)\lesssim a_{t_0}(x)^{1/2}.
$$
So for any $x\in E$,
$$
B_2(t,x)\lesssim \int_{t_0}^te^{-(\lambda_1-2\lambda_{\gamma(f)})s}\,ds a_{t_0}(x)^{1/2}\lesssim a_{t_0}(x)^{1/2}.
$$
Thus by the dominated convergence theorem, we get that for any $x\in E$,
\begin{equation*}
  \lim_{t\to\infty}B_2(t,x)
  =\int_{t_0}^\infty e^{-(\lambda_1-2\lambda_{\gamma(f)})s}e^{\lambda_1s}T_s(Af^2_1)(x)\,ds.
\end{equation*}
By \eqref{1.36}, we easily get that for any $x\in E$,
$$
B_3(t,x)\lesssim e^{-(\lambda_1-2\lambda_{\gamma(f)})t}a_{t_0}(x)^{1/2}\to 0
$$
as $t\to\infty$.
Thus, the proof of (3) is now complete. \hfill$\Box$

\section{Proofs of the Main Results}

In this section, we will prove the main results of this paper.
When referring to individuals in $X$ we will use the classical Ulam-Harris notation
so that every individual in
$X$ has a unique label, see \cite{HH}.
Although the Ulam-Harris labelling
of individuals is rich enough to encode genealogical order, the only feature we really need
of the Ulam-Harris notation is that individuals are uniquely identifiable amongst ${\cal T}$, the set of
labels of individuals realized in $X$.
For each individual
$u\in{\cal T}$ we shall write
$b_u$ and $d_u$ for its birth and death times respectively and
$\{z_u(r): r\in [b_u,d_u]\}$ for its spatial trajectory.
Define
$${\cal L}_t=\{ u\in {\cal T}, b_u\le t<d_u\},\quad t\ge 0.$$
Thus, $X_{s+t}$ has the following decomposition:
\begin{equation}\label{3.22}
  X_{s+t}=\sum_{u\in \cL_t}X^{u,t}_s,
\end{equation}
where given $\mathcal{F}_t$, $X^{u,t}_s$, $u\in \cL_t$, are independent and $X^{u,t}_s$ has the same law as $X_s$ under $\P_{\delta_{z_u(t)}}$.

\subsection{The large branching rate case: $\lambda_1>2\lambda_{\gamma(f)}$}

\begin{lemma}\label{lem:1.2}
If $\lambda_1>2\lambda_k$, then, for any $\nu\in {\cal M}_a(E)$, $H_t^{k,j}$ is a martingale under $\P_{\nu}$.
Moreover, the limit
\begin{equation}\label{1.12}
H_\infty^{k,j}:=\lim_{t\to\infty}H_t^{k,j}
\end{equation}
exists $\P_{\nu}$-a.s. and in $L^2(\P_{\nu})$.

\end{lemma}

\textbf{Proof:}
By the branching property, it suffices to prove the lemma for $\nu=\delta_x$ for $x\in E$.
Since $\phi_j^{(k)}(x)$ is an eigenfunction corresponding to $-\lambda_k$,
we have, for any $(t, x)\in (0, \infty)\times E$, $\P_{\delta_x} H_t^{k,j}=\phi_j^{(k)}(x)$.
Thus, by the Markov property, we get that, for any $x\in E$, $H_t^{k,j}$ is a martingale under $\P_{\delta_x}$.
Using \eqref{1.60}, we have that for any $x\in E$,
$$
\sup_{t>3t_0}\P_{\delta_x}(H_t^{k,j})^2\lesssim a_{t_0}(x)^{1/2}<\infty,
$$
from which the convergence asserted in the lemma follow easily.
\hfill$\Box$

Now we present the proof of Theorem \ref{The:1.2}.

\textbf{Proof of Theorem \ref{The:1.2}:}\quad
By the branching property, it suffices to prove the lemma for $\nu=\delta_x$ for $x\in E$.
  Define $M_t:=e^{\lambda_{\gamma(f)}t}\langle \widetilde{f}, X_t\rangle$.
It follows from the definition of $\widetilde{f}$ that $\gamma(\widetilde{f})\geq\gamma(f)+1$.
 From Lemma~\ref{lem:2.2}, we have the following:
\begin{description}
\item{(i)}
If $\lambda_1>2\lambda_{\gamma(\widetilde{f})}$, then for any $x\in E$,
\begin{eqnarray}
  \lim _{t\to\infty}e^{2\lambda_{\gamma(\tilde{f})}t}\P_{\delta_x}\langle \widetilde{f}, X_t\rangle^2
\end{eqnarray}
exists, thus we have
\begin{eqnarray*}
  \P_{\delta_x}M_t^2 &=& e^{-2(\lambda_{\gamma(\tilde f)}-\lambda_{\gamma(f)})t}e^{2\lambda_{\gamma(\tilde f)}t}
  \P_{\delta_x}\langle \widetilde{f}, X_t\rangle^2\\
   &=& O(e^{-2(\lambda_{\gamma(\tilde f)}-\lambda_{\gamma(f)})t})\to 0, \quad \mbox{ as } t\to\infty.
\end{eqnarray*}
\item{(ii)}
If $\lambda_1=2\lambda_{\gamma(\widetilde{f})}$,
then, for any $x\in E$, $\lim _{t\to\infty}t^{-1}e^{\lambda_1 t}
\P_{\delta_x}\langle \widetilde{f}, X_t\rangle^2$ exists.
Thus we have for any $x\in E$,
\begin{eqnarray*}
   \P_{\delta_x}M_t^2 &=& te^{-2(\lambda_{\gamma(\tilde f)}-\lambda_{\gamma(f)})t}(t^{-1}e^{\lambda_1 t}
   \P_{\delta_x}\langle \widetilde{f}, X_t\rangle^2)\\
   &=& O(te^{-2(\lambda_{\gamma(\tilde f)}-\lambda_{\gamma(f)})t})\to 0, \quad \mbox{ as } t\to\infty.
\end{eqnarray*}
\item{(iii)}
If $2\lambda_{\gamma(\widetilde{f})}>\lambda_1>2\lambda_{\gamma(f)}$,
then by Lemma \ref{lem:2.2}(2), for any $x\in E$,
$\lim _{t\to\infty}e^{\lambda_1 t}\P_{\delta_x}\langle \widetilde{f}, X_t\rangle^2$ exists.
Thus we have for any $x\in E$,
\begin{eqnarray*}
   \P_{\delta_x}M_t^2 &=& e^{-(\lambda_1-2\lambda_{\gamma(f)})t}(e^{\lambda_1 t} P_{\delta_x}
   \langle \widetilde{f}, X_t\rangle^2)\\
   &=& O(e^{-(\lambda_1-2\lambda_{\gamma(f)})t})\to 0, \quad \mbox{ as } t\to\infty.
\end{eqnarray*}
\end{description}
Combining the three cases above, we get that, for any $x\in E$,
 $\lim_{t\to\infty}M_t = 0$ in  $L^2(\P_{\delta_x})$.
Now using Lemma~\ref{lem:1.2}, we easily get the convergence in Theorem \ref{The:1.2}.\hfill$\Box$

\subsection{The small branching rate case: $\lambda_1<2\lambda_{\gamma(f)}$}

First, we recall some properties of weak convergence.
For $f:\mathbb{R}^d\to\mathbb{R}$, let $\|f\|_L:=\sup_{x\ne y}|f(x)-f(y)|/\|x-y\|$ and
  $\|f\|_{BL}:=\|f\|_{\infty}+\|f\|_L$. For any distributions $\nu_1$ and $\nu_2$ on $\mathbb{R}^d$, define
\begin{equation*}
  \beta(\nu_1,\nu_2):=\sup\left\{\left|\int f\,d\nu_1-\int f\,d\nu_2\right|~:~\|f\|_{BL}\leq1\right\}.
\end{equation*}
Then $\beta$ is a metric. By \cite[Theorem 11.3.3]{Dudley}, the topology generated by this metric
is equivalent to the weak convergence topology.
 From the definition, we can easily see that, if $\nu_1$ and $\nu_2$ are the distributions of
 two $\mathbb{R}^d$-valued random variables $X$ and $Y$ respectively, then
\begin{equation}\label{5.20}
  \beta(\nu_1,\nu_2)\leq \E\|X-Y\|\leq\sqrt{ \E\|X-Y\|^2}.
\end{equation}

\textbf{Proof of Theorem \ref{The:1.3}:}\quad
 We define an ${\mathbb R}^2$-valued random variable $U_1(t)$ by
\begin{equation}\label{2.4}
   U_1(t):=\left(e^{\lambda_1 t}\langle \phi_1,X_t\rangle, e^{\lambda_1t/2}\langle f, X_t\rangle\right).
\end{equation}
To get the conclusion of Theorem \ref{The:1.3}, it suffices to show that,
for any nonzero $\nu\in {\cal M}_a(E)$, under $\P_{\nu}$,
\begin{equation}\label{2.5a}
   U_1(t)\stackrel{d}{\to}\left(W_\infty, \sqrt{W_\infty}G_1(f)\right),
\end{equation}
where $G_1(f)\sim\mathcal{N}(0,\sigma_f^2)$ is independent of $W_\infty$.
To show the above, it suffices to show that,
for any $x\in E$, under $\P_{\delta_x}$,
\begin{equation}\label{2.5}
   U_1(t)\stackrel{d}{\to}\left(W_\infty, \sqrt{W_\infty}G_1(f)\right),
\end{equation}
where $G_1(f)\sim\mathcal{N}(0,\sigma_f^2)$ is independent of $W_\infty$.
In fact, if
$\nu=\sum_{j=1}^n\delta_{x_j}, n=1, 2, \dots,
\{x_j; j=1,\cdots, n\}\subset E$, then
$$
X_t=\sum_{j=1}^nX_t^j
$$
where $X^j_t$ is a branching Markov process starting from
$\delta_{x_j}, j=1, \dots, n$, and $X^j, j=1,\cdots, n$, are independent.
If \eqref{2.5} is valid, we put
$W^j_\infty:=\lim_{t\to\infty}e^{\lambda_1t}\langle \phi_1, X_t^j\rangle$.
Then we easily get under $\P_{\nu}$, $W_\infty=\sum_{j=1}^n W^j_\infty$.
For $\lambda_1<2\lambda_{\gamma(f)}$,
\begin{eqnarray*}
  &&\P_{\nu}\exp\left\{i\theta_1e^{\lambda_1 t}\langle \phi_1, X_t\rangle+i\theta_2e^{(\lambda_1/2) t}\langle f, X_t\rangle\right) \\
  &=&\prod_{j=1}^n  \P_{\nu}\exp\left\{i\theta_1e^{\lambda_1 t}\langle \phi_1, X_t^j\rangle+i\theta_2e^{(\lambda_1/2) t}\langle f, X_t^j\rangle\right) \\
  &\to& \prod_{j=1}^n\P_{\nu}\exp\left\{i\theta_1W^j_\infty-\frac{1}{2}\theta_2^2\sigma_f^2W^j_\infty\right) \\
  &=& \P_{\nu}\exp\left\{i\theta_1W_\infty-\frac{1}{2}\theta_2^2\sigma_f^2W_\infty\right),
\end{eqnarray*}
which implies that \eqref{2.5a} is valid.

Now we show that \eqref{2.5} is valid.
Let $s, t >3t_0$ and write
\begin{equation*}
  U_1(s+t)=\left(e^{\lambda_1 (s+t)}\langle \phi_1,X_{t+s}\rangle, e^{(\lambda_1/2)(s+t)}\langle f,X_{s+t}\rangle\right).
\end{equation*}
Recall the decomposition of $X_{s+t}$ in \eqref{3.22}.
Define
\begin{equation}\label{e:new}
Y_s^{u,t}:=e^{\lambda_1 s/2}\langle f,X_{s}^{u,t}\rangle\quad\mbox{and}\quad y_s^{u,t}:=\P_{\delta_x}(Y^{u,t}_s|\mathcal{F}_t).
\end{equation}
Given $\mathcal{F}_t$, $Y_s^{u,t}$ has the same law as $Y_s:=e^{\lambda_1 s/2}\langle f,X_{s}\rangle$ under $\P_{\delta_{z_u(t)}}$.
Then we have
\begin{eqnarray}\label{2.14}
  &&e^{(\lambda_1/2) (s+t)}\langle f,X_{s+t}\rangle= e^{(\lambda_1/2) t}\sum_{u\in\cL_t}Y_s^{u,t}\nonumber\\
  &=& e^{(\lambda_1/2) t}\sum_{u\in\cL_t}(Y_s^{u,t}-y^{u,t}_s)
      + e^{(\lambda_1/2) (t+s)}\P_{\delta_x}(\langle f,X_{s+t}\rangle|\mathcal{F}_t)\nonumber\\
  &=:& J_1(s,t)+J_2(s,t).
\end{eqnarray}

We first consider $J_{2}(s,t)$.
By the Markov property, we have
\begin{equation*}
  J_{2}(s,t)=e^{(\lambda_1/2)(s+t)}\langle T_sf, X_t\rangle.
\end{equation*}
Thus, by \eqref{1.13} and \eqref{1.36}, we have for any $x\in E$,
\begin{eqnarray*}
  \P_{\delta_x}(J_{2}(s,t)^2) &=& e^{\lambda_1(s+t)}\int_0^t T_{t-u}[A(T_{s+u}f)^2](x)\,du
  + e^{\lambda_1(s+t)}T_t[(T_sf)^2](x)\\
  &=&e^{(\lambda_1-2\lambda_{\gamma(f)})s}\int_0^t e^{(\lambda_1-2\lambda_{\gamma(f)})u}e^{\lambda_1(t-u)}T_{t-u}[A(e^{\lambda_{\gamma(f)}(s+u)}T_{s+u}f)^2](x)\,du\\
  &&+ e^{(\lambda_1-2\lambda_{\gamma(f)})s}e^{\lambda_1t}T_t[(e^{\lambda_{\gamma(f)}s}T_sf)^2](x)\\
  &\lesssim& e^{(\lambda_1-2\lambda_{\gamma(f)})s}
  \left(\int_0^t e^{(\lambda_1-2\lambda_{\gamma(f)})u}e^{\lambda_1(t-u)}T_{t-u}[a_{2t_0}](x)\,du+a_{t_0}(x)^{1/2}\right)
\end{eqnarray*}
and
\begin{eqnarray*}
  &&\int_0^t e^{(\lambda_1-2\lambda_{\gamma(f)})u}e^{\lambda_1(t-u)}T_{t-u}(a_{2t_0})(x)\,du\\
  &=&\left(\int_0^{t-t_0}+\int_{t-t_0}^t\right)
  e^{(\lambda_1-2\lambda_{\gamma(f)})u}e^{\lambda_1(t-u)}T_{t-u}(a_{2t_0})(x)\,du\\
  &\lesssim& \int_0^{t-t_0}e^{(\lambda_1-2\lambda_{\gamma(f)})u}\,du a_{2t_0}(x)^{1/2}
  +\int_0^{t_0}e^{(\lambda_1-2\lambda_{\gamma(f)})(t-u)}e^{\lambda_1u}T_{u}(a_{2t_0})(x)\,du\\
  &\lesssim& a_{t_0}(x)^{1/2}+\int_0^{t_0}T_{u}(a_{2t_0})(x)\,du\lesssim a_{t_0}(x)^{1/2}.
\end{eqnarray*}
Thus for any $x\in E$,
\begin{equation}\label{2.15}
  \limsup_{t\to\infty}\P_{\delta_x}(J_{2}(s,t)^2)\lesssim e^{(\lambda_1-2\lambda_{\gamma(f)})s}a_{t_0}(x)^{1/2}.
\end{equation}

Next we consider $J_1(s,t)$.
We define an ${\mathbb R}^2$-valued random variable $U_2(s,t)$ by
\begin{eqnarray*}
U_2(s,t):=\left(e^{\lambda_1 t}\langle \phi_1, X_t\rangle, J_1(s,t) \right).
\end{eqnarray*}
Let $V_s(x):=Var_{\delta_x}Y_s$. We claim that, for any $x\in E$, under $\P_{\delta_x}$,
\begin{equation}\label{2.1}
  U_2(s,t)\stackrel{d}{\to}\left(W_\infty, \sqrt{W_\infty}G_1(s)\right), \quad \mbox{ as } t\to\infty,
\end{equation}
where $G_1(s)\sim\mathcal{N}(0,\sigma^2_f(s))$ is independent of $W_\infty$ and $\sigma^2_f(s)=\langle V_s,\phi_1\rangle$.
Denote the characteristic function of $U_2(s,t)$ under $\P_{\delta_x}$ by
$\kappa(\theta_1,\theta_2,s,t)$:
\begin{eqnarray}
  \kappa(\theta_1,\theta_2,s,t)
  &=&\P_{\delta_x}\left(\exp\left\{i\theta_1e^{\lambda_1 t}\langle \phi_1, X_t\rangle+i\theta_2e^{(\lambda_1/2) t}\sum_{u\in\cL_t}(Y_s^{u,t}-y_s^{u,t})\right\}\right)\nonumber\\
  &=& \P_{\delta_x}\left(\exp\{i\theta_1e^{\lambda_1 t}\langle \phi_1, X_t\rangle\}\prod_{u\in\cL_t}h_s(z_u(t),e^{(\lambda_1/2) t}\theta_2)\right),\label{2.3}
\end{eqnarray}
where
$$
h_s(x,\theta)=\P_{\delta_x}e^{i\theta(Y_s-\P_{\delta_x}Y_s)}.
$$

Let $t_k,m_k\to\infty$, as $k\to\infty$, and $a_{k,j}\in E$, $j=1,2,\cdots m_k$. Now we consider
\begin{equation}\label{2.16}
  S_k:=e^{\lambda_1t_k/2}\sum_{j=1}^{m_k}(Y_{k,j}-y_{k,j}),
\end{equation}
where $Y_{k,j}$ has the same law as $Y_s$ under $\P_{\delta_{a_{k,j}}}$ and $y_{k,j}=\P_{\delta_{a_{k,j}}}Y_s$. Further, $Y_{k,j},j=1,2,\dots$ are independent. Suppose the Lindeberg conditions hold:
\begin{description}
\item{(i)} as $k\to\infty$,
    $$e^{\lambda_1t_k}\sum_{j=1}^{m_k}\E(Y_{k,j}-y_{k,j})^2=e^{\lambda_1t_k}\sum_{j=1}^{m_k}V_s(a_{k,j})\to\sigma^2;$$
  \item{(ii)} for any $\epsilon>0$,
   \begin{eqnarray*}
     && e^{\lambda_1t_k}\sum_{j=1}^{m_k}\E\left(|Y_{k,j}-y_{k,j}|^2,|Y_{k,j}-y_{k,j}|>\epsilon e^{-\lambda_1 t_k/2}\right) \\
     &=& e^{\lambda_1t_k}\sum_{j=1}^{m_k}\P_{\delta_{a_{k,j}}}\left(|Y_s-y_s|^2,|Y_s-y_s|>\epsilon e^{-\lambda_1 t_k/2}\right)\to 0,\quad \mbox{as }k\to\infty.
  \end{eqnarray*}
\end{description}
Then using the Lindeberg-Feller theorem,
we have $S_k\stackrel{d}{\to}\mathcal{N}(0,\sigma^2)$, which implies
\begin{equation}\label{2.17}
  \prod_{j=1}^{m_k}h_s(a_{k,j},e^{\lambda_1t_k/2}\theta)\to e^{-\frac{1}{2}\sigma^2\theta^2}.
\end{equation}
By Lemma \ref{lem:2.2}(1),
$V_s\in L^2(E,\mu)\cap L^4(E,\mu)$.
So using Remark \ref{rem:large}, we have
\begin{equation}\label{2.18}
  e^{\lambda_1t}\sum_{u\in \cL_t}V_s(z_u(t))=e^{\lambda_1t}\langle V_s,X_t\rangle\to\langle V_s,\phi_1\rangle W_\infty,\quad \mbox{  in probability}, \mbox{as }t\to\infty.
\end{equation}
Let $g(x,s,t)=\P_{\delta_x}\left(|Y_s-y_s|^2,|Y_s-y_s|>\epsilon e^{-\lambda_1 t/2}\right).$
We note that $g(x,s,t)\downarrow0$ as  $t\uparrow\infty$ and $g(x,s,t)\leq V_s(x)$ for any $x\in E$.
Thus by Lemma \ref{lem:rsnew} we have for any $x\in E$,
\begin{equation*}
 e^{\lambda_1t}\P_{\delta_x}\langle g(\cdot,s,t),X_t\rangle=e^{\lambda_1t}T_t(g(\cdot,s,t))(x)\to 0,
 \quad \mbox{as}\quad t\to\infty,
\end{equation*}
which implies
\begin{equation}\label{2.19}
  e^{\lambda_1t}\sum_{u\in \cL_t}\P_{\delta_{z_u(t)}}\left(|Y_s-y_s|^2,|Y_s-y_s|>\epsilon e^{-\lambda_1 t/2}\right)\to 0,\quad \mbox{as}\quad t\to\infty,
\end{equation}
in $\P_{\delta_x}$-probability.
Therefore, for any sequence $s_k\to\infty$, there exists a subsequence $s_k'$ such that,
if we let $t_k=s_k'$, $m_k=|X_{s_k'}|$ and
$\{a_{k,j},j=1,2\cdots m_k\}=\{z_u(s_k'),u\in \cL_{s_k'}\}$,
then the Lindeberg conditions hold $\P_{\delta_x}$-a.s. for any $x\in E$, which implies
\begin{equation}\label{2.25}
  \lim_{k\to\infty}\prod_{u\in\cL_{s_k'}}h_s(z_u(s_k'),e^{(\lambda_1/2) s_k'}\theta_2)
  =\exp\left\{-\frac{1}{2}\theta_2^2\langle V_s,\phi_1\rangle W_\infty\right\},\quad\P_{\delta_x}\mbox{-a.s.}
\end{equation}
Consequently, we have
\begin{equation}\label{2.20}
  \lim_{t\to\infty}\prod_{u\in\cL_t}h_s(z_u(t),e^{(\lambda_1/2) t}\theta_2)
  =\exp\left\{-\frac{1}{2}\theta_2^2\langle V_s,\phi_1\rangle W_\infty\right\},\quad\mbox{in probability}.
\end{equation}
Hence by the dominated convergence theorem, we get
\begin{equation}\label{2.21}
 \lim_{t\to\infty}\kappa(\theta_1,\theta_2,s,t)= \P_{\delta_x}\exp\left\{i\theta_1W_\infty\right\}
  \exp\left\{-\frac{1}{2}\theta_2^2\langle V_s,\phi_1\rangle W_\infty\right\},
\end{equation}
which implies our claim \eqref{2.1}.
Thus, we easily get that, for any $x\in E$, under $\P_{\delta_x}$,
\begin{eqnarray*}
U_3(s,t):=\left(e^{\lambda_1 (t+s)}\langle \phi_1, X_{t+s}\rangle,J_1(s,t) \right)
\stackrel{d}{\to}\left(W_\infty, \sqrt{W_\infty}G_1(s)\right), \quad \mbox{ as } t\to\infty.
\end{eqnarray*}
By \eqref{limit-mean} and \eqref{2.2}, we have
$$
\lim_{s\to\infty}\langle V_s,\phi_1\rangle=\sigma^2_f.
$$
Let $G_1(f)$ be a  $\mathcal{N}(0,\sigma_f^2)$ random variable independent of $W_\infty$.
Then
\begin{equation}\label{2.22}
\lim_{s\to\infty}\beta(G_1(s),G_1(f))=0.
\end{equation}
Let $\mathcal{L}(s+t)$ and $\widetilde{\mathcal{L}}(s,t)$ be the distributions of $U_1(s+t)$ and $U_3(s,t)$
respectively, and let $\mathcal{L}(s)$ and $\mathcal{L}$ be the distributions of $(W_\infty, \sqrt{W_\infty}G_1(s))$
and $(W_\infty, \sqrt{W_\infty}G_1(f))$ respectively. Then, using \eqref{5.20}, we have
\begin{eqnarray}\label{2.23}
  \limsup_{t\to\infty}\beta(\mathcal{L}(s+t),\mathcal{L})&\leq&
  \limsup_{t\to\infty}[\beta(\mathcal{L}(s+t),\widetilde{\mathcal{L}}(s,t))+\beta(\widetilde{\mathcal{L}}(s,t),\mathcal{L}(s))
  +\beta(\mathcal{L}(s),\mathcal{L})]\nonumber\\
 &\leq &\limsup_{t\to\infty}\P_{\delta_x}(J_2(s,t)^2)^{1/2}+0+\beta(\mathcal{L}(s),\mathcal{L}).
\end{eqnarray}
Using this and the definition of $\limsup_{t\to\infty}$, we easily get that
$$
\limsup_{t\to\infty}\beta (\mathcal{L}(t),\mathcal{L})=\limsup_{t\to\infty}\beta (\mathcal{L}(s+t),\mathcal{L})
\le \limsup_{t\to\infty}(\P_{\delta_x}J_2(s,t)^2)^{1/2}+\beta(\mathcal{L}(s),\mathcal{L}).
$$
Letting $s\to\infty$, we get $ \limsup_{t\to\infty}\beta(\mathcal{L}(t),\mathcal{L})=0$.
The proof is now complete. \hfill$\Box$

\subsection{Proof of Theorem \ref{The:2.1}}

In this subsection we consider the case: $\lambda_1>2\lambda_{\gamma(f)}$ and $f_{(c)}=0$.
By Lemma \ref{lem:1.2}, we have for $2\lambda_k<\lambda_1$,
\begin{equation}\label{3.1}
  H_\infty^{k,j}=\lim_{s\to\infty}e^{\lambda_k(t+s)}\langle \phi_j^{(k)},X_{t+s}\rangle
  =e^{\lambda_kt}\lim_{s\to\infty}\sum_{u\in\cL_t}\langle \phi_j^{(k)},X_{s}^{u,t}\rangle e^{\lambda_ks}.
\end{equation}
Since $X_{s}^{u,t}$ has the same law as $X_s$ under $\P_{\delta_{z_u(t)}}$,
$H^{u,t,k,j}_\infty=\lim_{s\to\infty}\langle \phi_j^{(k)},X_{s}^{u,t}\rangle e^{\lambda_ks}$ exists
and has the same law as $H^{k,j}_\infty$ under $\P_{\delta_{z_u(t)}}$.
Thus
\begin{equation}\label{3.2}
  H^{k,j}_\infty=e^{\lambda_kt}\sum_{u\in\cL_t}H^{u,t,k,j}_\infty.
\end{equation}
Recall from \eqref{1.27} that
$$
H_\infty:=\sum_{2\lambda_k<\lambda_1}\sum_{j=1}^{n_k}a_j^kH^{k,j}_\infty.
$$
Denote
$$
H_\infty^{u,t}:=\sum_{2\lambda_k<\lambda_1}\sum_{j=1}^{n_k}a_j^kH^{u,t,k,j}_\infty.
$$
It is easy to see that given $\mathcal{F}_t$, $H_\infty^{u,t}$ has the same law as $H_\infty$ under $\P_{\delta_{z_u(t)}}$.
By Lemma \ref{lem:1.2}, we have that for any $x\in E$,
\begin{equation*}
  \sum_{2\lambda_k<\lambda_1}\sum_{j=1}^{n_k}e^{\lambda_kt}a_j^k\langle
  \phi_j^{(k)},X_{t}\rangle\to H_\infty, \quad\mbox{ in }L^2(\P_{\delta_x}).
\end{equation*}
It follows that
\begin{equation}\label{3.3}
  \P_{\delta_x}H_\infty=f_{(s)}(x), \qquad x\in E
\end{equation}
and by \eqref{1.13}, we have that for any $x\in E$,
\begin{equation}\label{3.4}
  \P_{\delta_x}(H_\infty)^2=\int_0^\infty T_s
    \left[A\left(\sum_{2\lambda_k<\lambda_1}\sum_{j=1}^{n_k}e^{\lambda_ks}a_j^k\phi_j^k\right)^2\right](x)\,ds.
\end{equation}

 \textbf{Proof of Theorem \ref{The:2.1}:}\quad
 By \eqref{3.2}, we have
 $$
 \sum_{\lambda_k<\lambda_1/2}e^{-\lambda_kt}\sum_{j=1}^{n_k}a_j^kH^{k,j}_\infty=\sum_{u\in\cL_t}H^{u,t}_\infty.
 $$
 Consider the $\mathbb{R}^2$-valued random variable $U_1(t)$:
 \begin{equation}\label{3.5}
   U_1(t):= \left(e^{\lambda_1 t}\langle \phi_1,X_t\rangle,~e^{(\lambda_1/2)t}\left(\langle
   f,X_t\rangle-\sum_{u\in{\cL}_t}H^{u,t}_\infty\right)\right).
 \end{equation}
Using an argument similar to that in the beginning of the proof of
Theorem \ref{The:1.3}, we can see that, to get the conclusion of Theorem \ref{The:2.1}, it suffices to show that
for any $x\in E$, under $\P_{\delta_x}$, as $t\to\infty$,
\begin{equation}\label{3.6}
 U_1(t)\stackrel{d}{\to}\left(W_\infty, \sqrt{W_\infty}G_3(f)\right),
\end{equation}
where
 $G_3(f)\sim \mathcal{N}(0,\sigma_{f_{(l)}}^2+\beta_{f}^2)$ is independent of $W_\infty$.
 Denote the characteristic function of $U_1(t)$  under  $\P_{\delta_x}$ by $\kappa_1(\theta_1,\theta_2,t)$ and let
 $h(x,\theta):=\P_{\delta_x}\exp\{i\theta (H_\infty-f_{(s)}(x))\}$.
 Then we have for any $x\in E$,
 \begin{eqnarray}\label{3.7}
  &&\kappa_1(\theta_1,\theta_2,t)\nonumber\\
  &=&\P_{\delta_x}\exp\left\{i\theta_1e^{\lambda_1 t}\langle \phi_1,X_t\rangle+i\theta_2e^{(\lambda_1/2)t}\left(\langle
   f,X_t\rangle-\sum_{u\in{\cL}_t}H^{u,t}_\infty\right)\right\}\nonumber\\
   &=&\P_{\delta_x}\exp\left\{i\theta_1e^{\lambda_1 t}\langle \phi_1,X_t\rangle+i\theta_2e^{(\lambda_1/2)t}\left(\langle
   f_{(l)},X_t\rangle-\sum_{u\in{\cL}_t}\left(H^{u,t}_\infty-f_{(s)}(z_u(t))\right)\right)\right\}\nonumber\\
   &=&\P_{\delta_x}\exp\left\{i\theta_1e^{\lambda_1 t}\langle \phi_1,X_t\rangle\right\}\exp\left\{i\theta_2e^{(\lambda_1/2)t}\langle f_{(l)},X_t\rangle\right\}\prod_{u\in{\cL}_t}
   h\left(z_u(t),-\theta_2e^{(\lambda_1/2)t}\right).
 \end{eqnarray}
Let $V(x)=Var_{\delta_x}H_\infty$. We claim that
\begin{description}
  \item{(i)}
  as $t\to\infty$,
  \begin{equation}\label{3.24}
    e^{\lambda_1t}\sum_{u\in{\cL}_t}\P_{\delta_x}(H^{u,t}_\infty-f_{(s)}(z_u(t)))^2
  =e^{\lambda_1t}\langle V,X_t\rangle\to\langle V,\phi_1\rangle W_\infty, \mbox{ in probability};
  \end{equation}
  \item{(ii)} for any $\epsilon>0$,
  as $t\to\infty$,
  \begin{equation}\label{3.10}
    e^{\lambda_1t}\sum_{u\in{\cL}_t}
    \P_{\delta_x}(|H^{u,t}_\infty-f_{(s)}(z_u(t))|^2,|H^{u,t}_\infty-f_{(s)}(z_u(t))|>\epsilon e^{-\lambda_1t/2})\to 0,
    \mbox{ in probability.}
  \end{equation}
\end{description}
Then using arguments similar to those in the proof Theorem \ref{The:1.3}, we have
\begin{equation}\label{3.11}
  \prod_{u\in{\cL}_t}h\left(z_u(t),-\theta_2e^{(\lambda_1/2)t}\right)
  \to \exp\left\{-\frac{1}{2}\theta_2^2\langle V,\phi_1\rangle W_\infty\right\},\mbox{ in probability.}
\end{equation}
Now we will prove the claims.

(i)
 By Remark \ref{rem:large}, we only need to show that $V(x)\in L^2(E,\mu)\cap L^4(E,\mu)$.
By \eqref{1.37}, we have that for any $x\in E$,
$$
\sum_{2\lambda_k<\lambda_1}\sum_{j=1}^{n_k}e^{\lambda_ks}|a_j^k||\phi_j^k(x)|\lesssim e^{\lambda_m s}a_{2t_0}(x)^{1/2},
$$
 where $m=\sup\{k:2\lambda_k<\lambda_1\}$.
 So by \eqref{3.4} and \eqref{1.36}, we have that for any $x\in E$,
 \begin{eqnarray*}
  \P_{\delta_x}(H_\infty)^2&\lesssim& \int_0^\infty e^{(2\lambda_m-\lambda_1)s}e^{\lambda_1s}T_s(a_{2t_0})(x)\,ds\\
    &=&\left(\int_0^{t_0}+\int_{t_0}^\infty\right) e^{(2\lambda_m-\lambda_1)s}e^{\lambda_1s}T_s(a_{2t_0})(x)\,ds\\
    &\lesssim& \int_0^{t_0}T_s(a_{2t_0})(x)\,ds+\int_{t_0}^\infty e^{(2\lambda_m-\lambda_1)s}\,ds~a_{t_0}(x)^{1/2}\\
    &\lesssim& a_{t_0}(x)^{1/2}\in L^2(E,\mu)\cap L^4(E,\mu).
\end{eqnarray*}
Thus $V(x)\in L^2(E,\mu)\cap L^4(E,\mu)$.

(ii)
Let
$$
g_t(x)=\P_{\delta_x}\left(|H_\infty-f_{(s)}(x)|^2,|H_\infty-f_{(s)}(x)|>\epsilon e^{-\lambda_1t/2}\right), \qquad (t, x)\in (0, \infty)\times E.
$$
Then for any $t>0$,
$$
e^{\lambda_1t}\sum_{u\in{\cL}_t}
\P_{\delta_x}\left(|H^{u,t}_\infty-f_{(s)}(z_u(t))|^2,|H^{u,t}_\infty-f_{(s)}(z_u(t))|>\epsilon e^{-\lambda_1t/2}\right)
=e^{\lambda_1t}\langle g_t,X_t\rangle.
$$
We easily see that $g_t(x)\downarrow0$ as  $t\uparrow\infty$ and $g_t(x)\leq V(x)$ for any $x\in E$.
So, by Lemma \ref{lem:rsnew}, we have that for any $x\in E$,
\begin{equation*}
  e^{\lambda_1t}\P_{\delta_x}\langle g_t,X_t\rangle=e^{\lambda_1t}T_t(g_t)(x)\to 0,\quad\mbox{ as }t\to\infty,
\end{equation*}
which implies \eqref{3.10}.

By \eqref{3.11} and the dominated convergence theorem, we get that as $t\to\infty$,
\begin{equation}\label{3.12}
  \left|\kappa_1(\theta_1,\theta_2,t)
  -\P_{\delta_x}\exp\left\{i\theta_1e^{\lambda_1 t}\langle \phi_1,X_t\rangle+i\theta_2e^{\lambda_1t/2}
  \langle   f_{(l)},X_t\rangle-\frac{1}{2}\theta_2^2\langle V,~\phi_1\rangle W_\infty\right\}\right|\to 0.
\end{equation}
Since $\lambda_1<2\lambda_{\gamma(f_{(l)})}$, by Theorem \ref{The:1.3}, we have that as $t\to\infty$,
\begin{equation}\label{3.13}
  \left(e^{\lambda_1 t}\langle \phi_1,X_t\rangle,e^{(\lambda_1/2)t}\langle f_{(l)},X_t\rangle\right)\stackrel{d}{\to}\left(W_\infty,\sqrt{W_\infty}G_1(f_{(l)})\right),
\end{equation}
where $G_1(f_{(l)})\sim\mathcal{N}(0,\sigma^2_{f_{(l)}})$ is independent of $W_\infty$.
 Therefore, for any $x\in E$, as $t\to\infty$,
 \begin{eqnarray}\label{3.14}
   &&\P_{\delta_x}\exp\left\{i\theta_1e^{\lambda_1 t}\langle \phi_1,X_t\rangle+i\theta_2e^{(\lambda_1/2)t}\langle f_{(l)},X_t\rangle\right\}\exp\left\{-\frac{1}{2}\theta_2^2\langle V,~\phi_1\rangle W_\infty\right\}\nonumber\\
   &\to& \P_{\delta_x}\left(\exp\{i\theta_1W_\infty\}\exp\{-\frac{1}{2}\theta_2^2(\sigma^2_{f_{(l)}}+\langle V,~\phi_1\rangle )W_\infty\}\right).
\end{eqnarray}
By \eqref{3.3} and \eqref{3.4}, we get
\begin{equation*}
   \langle V,~\phi_1\rangle =\int_0^\infty e^{-\lambda_1 s}
    \left\langle A(\sum_{2\lambda_k<\lambda_1}\sum_{j=1}^{n_k}e^{\lambda_ks}a_j^k\phi_j^k)^2,\phi_1\right\rangle\,ds
    -\langle (f_{(s)})^2,\phi_1\rangle.
\end{equation*}
The proof is now complete.\hfill$\Box$

\subsection{The critical branching rate case: $\lambda_1=2\lambda_{\gamma(f)}$}

To prove Theorem \ref{The:1.4}, we need the following lemma.
\begin{lemma}\label{lem:5.5}
 Assume $f=\sum_{j=1}^{n_k}b_j^k\phi_j^{(k)}$, where $b_j^k\in \mathbb{R}$ and $\lambda_1=2\lambda_k$.
Define
$$
  S_tf(x):=t^{-1/2}e^{(\lambda_1/2) t}(\langle f,X_t\rangle-{T}_tf(x)), \qquad (t, x)\in (0, \infty)\times E.
$$
Then for any $c>0$, $\delta>0$ and $x\in E$, we have
\begin{equation}\label{4.5}
 \lim_{t\to\infty}\mathbb{P}_{\delta_x}\left(|S_tf(x)|^2;|S_tf(x)|>ce^{\delta t }\right)=0.
\end{equation}
\end{lemma}

\textbf{Proof:}\quad
We write $t=[t]+\epsilon_t$, where $[t]$ is the integer part of $t$.
Let
 $$
 F(t, x):=\P_{\delta_x}\left(|S_tf(x)|^2;|S_tf(x)|>ce^{\delta t }\right),  \qquad (t, x)\in (0, \infty)\times E.
 $$
By the definition of $f$, we easily get $T_uf(x)=e^{-\lambda_1u/2}f(x)$.
So we get that for any $(t, x)\in (0, \infty)\times E$,
\begin{eqnarray}\label{4.8}
S_{t+1}f(x)&=&\left(\frac{1}{t+1}\right)^{1/2}e^{(\lambda_1/2) (t+1)}\left(\langle f, X_{t+1}\rangle-
  \langle e^{-\lambda_1/2}f, X_t\rangle\right)\nonumber\\
  &&+\left(\frac{1}{t+1}\right)^{1/2}e^{(\lambda_1/2)t}\left(\langle f, X_t\rangle-T_{t}f(x)\right)\nonumber\\
   &=&\left(\frac{1}{t+1}\right)^{1/2}R(t,f)+\left(\frac{t}{t+1}\right)^{1/2}S_tf(x),
\end{eqnarray}
where
$$
R(t,f):=e^{(\lambda_1/2) (t+1)}\left(\langle f,X_{t+1}\rangle-\langle T_1f, X_t\rangle\right).
$$
Thus we have that for any $(t, x)\in (0, \infty)\times E$,
 \begin{eqnarray*}
 &&F(t+1, x)\\
 &\leq& \P_{\delta_x}\left(|S_{t+1}f(x)|^2;|S_tf(x)|>ce^{\delta t }\right)
 +\P_{\delta_x}\left(|S_{t+1}f(x)|^2;|S_tf(x)|\leq ce^{\delta t},|S_{t+1}f(x)|>ce^{\delta (t+1)}\right)\\
 &=:&M_1(t, x)+M_2(t, x).
 \end{eqnarray*}
 Put
 \begin{eqnarray*}
    A_1(t,x) &=& \{|S_tf(x)|>ce^{\delta t }\},\\
    A_2(t,x) &=& \{|S_tf(x)|\leq ce^{\delta t }, |S_{t+1}f(x)|> ce^{\delta(t+1)}\}.
 \end{eqnarray*}
 Since
$A_1(t, x)\in\mathcal{F}_t$ and $\P_{\delta_x}(R(t,f) |\mathcal{F}_t)$=0 for any $(t, x)\in (0, \infty)\times E$,
we have by \eqref{4.8} that
 \begin{eqnarray*}
 M_1(t, x)&=&\frac{1}{t+1}\P_{\delta_x} \left(|R(t,f)|^2;A_1(t,x)\right) +\frac{t}{t+1}F(t, x),
 \end{eqnarray*}
 and
 \begin{eqnarray*}
   M_2(t, x)&\leq&\frac{2}{t+1}\P_{\delta_x}\left(|R(t,f)|^2;A_2(t,x)\right)
   +\frac{2t}{t+1}\P_{\delta_x}\left(|S_tf(x)|^2;A_2(t, x)\right).
 \end{eqnarray*}
 Thus we have that  for any $(t, x)\in (0, \infty)\times E$,
 \begin{eqnarray}\label{4.9}
    F(t+1, x)\leq\frac{t}{t+1}F(t, x)+\frac{1}{t+1}(F_1(t, x)+F_2(t, x)),
 \end{eqnarray}
 where
 \begin{eqnarray*}
   F_1(t, x) &=& 2\P_{\delta_x}\left(|R(t,f)|^2;A_1(t, x)\cup A_2(t, x)\right),\\
   F_2(t, x) &=& 2t\P_{\delta_x}\left(|S_tf(x)|^2;A_2(t, x)\right).
 \end{eqnarray*}
 Choose an integer $k_0>3t_0$. Iterating \eqref{4.9}, we get for $t$ large enough
 \begin{eqnarray}\label{4.10}
    F(t+1, x)&\leq&\frac{1}{t+1}\sum_{m=k_0}^{[t]}\left(F_1(m+\epsilon_t,x)+F_2(m+\epsilon_t,x)\right)
    +\frac{k_0+\epsilon_t}{t+1}F(k_0+\epsilon_t,x)\nonumber\\
    &:=&L_1(t, x)+L_2(t, x)+\frac{k_0+\epsilon_t}{t+1}F(k_0+\epsilon_t,x).
 \end{eqnarray}

First, we will consider $L_1(t, x)$.
By \eqref{1.36}, we have that for any $x\in E$ and $s\ge k_0$,
\begin{eqnarray}\label{G1}
F_1(s, x) \le 2\P_{\delta_x}\left(|R(s,f)|^2\right)
=2e^{\lambda_1(s+1)}T_s(Var_{\delta_\cdot}\langle f,X_1\rangle)(x)
\le Ca_{t_0}(x)^{1/2},
\end{eqnarray}
where $C$ is a constant.
We claim that for any $x\in E$,
\begin{equation}\label{Mto0}
 F_1(t, x)\to 0,\quad \mbox{ as } t\to\infty.
\end{equation}
Then, for any $\epsilon>0$ and $x\in E$, there exists $K\in \mathbb{N}$
such that $s\ge K$ implies $F_1(s, x)<\epsilon$.
So, by \eqref{G1}, we get that for any $x\in E$ and $t$ large enough,
\begin{eqnarray*}
  L_1(t,x)&=& \frac{1}{t+1}\sum_{m=k_0}^{K-1}F_1(m+\epsilon_t,x)+\frac{1}{t+1}\sum_{m=K}^{[t]}F_1(m+\epsilon_t,x)
  \le \frac{CK}{t+1}a_{t_0}(x)^{1/2}+\epsilon.
\end{eqnarray*}
 Thus $\limsup_{t\to\infty}L_1(t,x)\le \epsilon$ for any $x\in E$,
which implies
\begin{equation}\label{L1}
  \lim_{t\to\infty}L_1(t, x)=0, \qquad x\in E.
\end{equation}

Now we prove the claim.
First, we will show that, for any $x\in E$, as $t\to\infty$,
\begin{equation}\label{Ato0}
   \P_{\delta_x}(A_1(t, x)\cup A_2(t, x))\to 0.
\end{equation}
By Chebyshev's inequality and \eqref{1.49}, we have that, for any $x\in E$, as $t\to\infty$,
\begin{eqnarray*}
  &&\P_{\delta_x}( A_1(t, x)) \leq c^{-2}e^{-2\delta t }\P_{\delta_x}|S_tf(x)|^2\to 0.
\end{eqnarray*}
It is easy to see that, under $\P_{\delta_x}$, for any $t>0$,
\begin{equation}\label{4.24}
A_2(t, x)\subset\left\{|R(t,f)|>ce^{\delta t}(e^{\delta}\sqrt{t+1}-\sqrt{t})\right\}.
\end{equation}
Similarly, by Chebyshev's inequality, we have that, for any $x\in E$,
\begin{equation*}
 \P_{\delta_x}(A_2(t, x))\leq c^{-2}e^{-2\delta t}(e^{\delta}\sqrt{t+1}-\sqrt{t})^{-2}\P_{\delta_x}|R(t,f)|^2.
\end{equation*}
By \eqref{1.25}, we get that, for any $x\in E$,
\begin{equation}\label{var:R}
  \P_{\delta_x}|R(t,f)|^2=e^{\lambda_1(t+1)}T_t(Var_{\delta_\cdot}\langle f,X_1\rangle)(x)
  \to e^{\lambda_1}\langle Var_{\delta_\cdot}\langle f,X_1\rangle,\phi_1\rangle\phi_1(x),
\end{equation}
which implies $\P_{\delta_x}(A_2(t, x))\to 0$ for any $x\in E$.

Using \eqref{2.14}, we have
\begin{eqnarray*}
  R(t,f)&=&e^{(\lambda_1/2) (t+1)}\left(\langle f,X_{t+1}\rangle-\langle T_1f,X_t\rangle\right)
  =e^{(\lambda_1/2) t}\sum_{u\in\cL_t}(Y_1^{u,t}-y^{u,t}_1),
\end{eqnarray*}
where $Y_1^{u,t}, y^{u,t}_1$ are defined in \eqref{e:new}.  From the proof of \eqref{2.1}, we see that
 \eqref{2.1} is also true when $\lambda_1=2\lambda_{\gamma(f)}$.
 Recall that $V_1(x)=e^{\lambda_1}Var_{\delta_x}\langle f,X_1\rangle$ for any $x\in E$.
So we have $R(t,f)\stackrel{d}{\to}\sqrt{W_\infty}G$,
where $G\sim\mathcal{N}(0, \langle V_1,\phi_1\rangle)$ is independent of $W_\infty$.
Then $\P_{\delta_x}(W_\infty G^2)=\langle V_1,\phi_1\rangle\phi_1(x)$ for any $x\in E$.

Let $\Psi_{M}(r)=r$ on $[0,M-1]$, $\Psi_M(r)=0$ on $[M,\infty]$, and let $\Psi_M$ be linear on $[M-1,M]$.
Therefore by \eqref{Ato0} and \eqref{var:R}, we have that for any $x\in E$,
\begin{eqnarray*}
\limsup_{t\to\infty}F_1(t, x) &\le& \limsup_{t\to\infty}2 \P_{\delta_x}(|R(t,f)|^2;|R(t,f)|^2>M)
   +2M\limsup_{t\to\infty}\P_{\delta_x}(A_1(t,x)\cup A_2(t,x)) \\
   &\le& 2 \limsup_{t\to\infty}\left(\P_{\delta_x}(|R(t,f)|^2)-\P_{\delta_x}(\Psi_M(|R(t,f)|^2))\right)\\
   &=&2 \left(\langle V_1,\phi_1\rangle\phi_1(x)-\P_{\delta_x}(\Psi_M(W_\infty G^2))\right).
\end{eqnarray*}
By the monotone convergence theorem, we have that for any $x\in E$,
$$\lim_{M\to\infty}\P_{\delta_x}(\Psi_M(W_\infty G^2))=\P_{\delta_x}(W_\infty G^2)=\langle V_1,\phi_1\rangle\phi_1(x),$$
which implies $F_1(t, x)\to 0$ for any $x\in E$.

Now we consider $L_2(t, x)$. We also claim that for any $x\in E$,
\begin{equation}\label{G2to0}
  F_2(t,x)\to 0, \quad\mbox{as}\quad t\to\infty.
\end{equation}
In fact, by \eqref{4.24}, we have that for any $x\in E$,
\begin{eqnarray*}
  F_2(t, x)&=& 2t\P_{\delta_x}\left(|S_tf(x)|^2;A_2(t,x)\right)\\
 &\leq& 2tce^{\delta t} \P_{\delta_x}\left(|S_tf(x)|;|R(t,f)|>ce^{\delta t}(e^{\delta}\sqrt{t+1}-\sqrt{t})\right)\\
    &\leq& 2c^{-1}te^{-\delta t}(e^{\delta}\sqrt{t+1}-\sqrt{t})^{-2}\P_{\delta_x}\left(|S_tf(x)|\cdot|R(t,f)|^2\right)\\
   &\lesssim& e^{-\delta t}e^{\lambda_1(t+1)}\P_{\delta_x}\left(|S_tf(x)|\langle Var_{\delta_\cdot}
   \langle f,X_1\rangle,X_t\rangle\right)\\
   &\lesssim& e^{-\delta t}\sqrt{\P_{\delta_x}|S_tf(x)|^2}\sqrt{e^{2\lambda_1t}\P_{\delta_x}
   \left(\langle Var_{\delta_\cdot}\langle f,X_1\rangle,X_t\rangle^2\right)}.
\end{eqnarray*}
By \eqref{1.49} and \eqref{1.23}, we get $F_2(t, x)\to 0$ for any $x\in E$ as $t\to \infty$.
Thus, for any $\epsilon>0$ and $x\in E$, there exists $K\in \mathbb{N}$ such that $s\ge K$ implies $F_2(s,x)<\epsilon$.
It is easy to see that,
$$
 \sup_{s\le K}F_2(s,x)\le \sup_{s\le K}2c^2se^{2\delta s}\le 2c^2Ke^{2\delta K}.
$$
Thus, we get
\begin{eqnarray*}
  L_2(t,x)=\frac{1}{t+1}\sum_{m=k_0}^{K-1}F_2(m+\epsilon_t,x)+\frac{1}{t+1}\sum_{m=K}^{[t]}F_2(m+\epsilon_t,x)
  \le \frac{2c^2K^2e^{2\delta K}}{t+1}+\epsilon.
\end{eqnarray*}
Thus $\limsup_{t\to\infty}L_2(t,x)\le \epsilon,$ which implies
\begin{equation}\label{L2}
  \lim_{t\to\infty}L_2(t,x)=0.
\end{equation}

To finish the proof,  we need to show that for any $x\in E$,
\begin{equation}\label{L3}
\lim_{t\to\infty}\frac{k_0+\epsilon_t}{t+1}F(k_0+\epsilon_t,x)=0.
\end{equation}
By \eqref{1.49}, we get that for any $x\in E$,
\begin{eqnarray*}
   \sup_{t>0}(k_0+\epsilon_t)F(k_0+\epsilon_t,x) \le
  (k_0+1)\sup_{s\ge k_0}\P_{\delta_x}(S_{s}f(x))^2<\infty,
\end{eqnarray*}
which implies \eqref{L3}.

Thus, we finish the proof. \hfill$\Box$

\smallskip

Now we are ready to prove Theorem \ref{The:1.4}.

\textbf{Proof of Theorem \ref{The:1.4}:}\quad
The proof is similar to that of Theorem \ref{The:1.3}.
 We define an ${\mathbb R}^2$-valued random variable by
\begin{equation*}
  U_1(t):=(e^{\lambda_1 t}\langle \phi_1, X_t\rangle, t^{-1/2}e^{(\lambda_1/2) t}\langle f_1,X_{t}\rangle).
\end{equation*}
 Since $\lambda_1=2\lambda_{\gamma(f)}$, $f=f_1+f_{(l)}$.
Using Theorem \ref{The:1.3} for $f_{(l)}$, we have
$$t^{-1/2}e^{(\lambda_1/2) t}\langle f_{(l)},X_{t}\rangle \stackrel{d}{\to} 0,\quad t\to \infty.$$
Thus, using an argument similar to that in the beginning of the proof of Theorem \ref{The:1.3},
to get conclusion of Theorem \ref{The:1.4}, we only need to show that, for any $x\in E$,
under $\P_{\delta_x}$, as $t\to\infty$,
\begin{equation}\label{4.1}
  U_1(t)\stackrel{d}{\to}\left(W_\infty, \sqrt{W_\infty}G_2(f)\right),
  \end{equation}
 where $G_2(f)\sim\mathcal{N}(0,\rho_f^2)$ is independent of $W_\infty$.
Let $t>3t_0$ and $n>2$. We write
\begin{equation*}
  U_1(nt)=\left(e^{\lambda_1 (nt)}\langle \phi_1,X_{nt}\rangle, (nt)^{-1/2}e^{(\lambda_1/2) (nt)}\langle f_1,X_{nt}\rangle\right).
\end{equation*}
Define
$$
Y_t^{u, n}:=((n-1)t)^{-1/2}e^{\lambda_1(n-1)t/2}\left\langle f_1,X_{(n-1)t}^{u,t}\right\rangle.
$$
Given $\mathcal{F}_t$, $Y_t^{u, n}$ has the same distribution as $Y^n_t:=((n-1)t)^{-1/2}e^{\lambda_1(n-1)t/2}\langle f_1,X_{(n-1)t}\rangle$
under $\P_{\delta_{Z_u(t)}}$. Since for $u>0$, $T_u f_1(x)=e^{-\lambda_1 u/2}f_1(x)$, we have
$$y_t^{u, n}:=\P_{\delta_x}(Y^{u, n}_t|\mathcal{F}_t)=((n-1)t)^{-1/2}f_1(z_u(t)).$$
Thus
\begin{eqnarray}\label{4.2}
     &&(nt)^{-1/2}e^{(\lambda_1/2)nt}\langle f_1,X_{nt}\rangle= \sqrt{\frac{n-1}{n}}e^{(\lambda_1/2) t}\sum_{u\in\cL_t}Y_t^{u, n}\nonumber\\
  &=& \sqrt{\frac{n-1}{n}}e^{(\lambda_1/2) t}\sum_{u\in\cL_t}(Y_t^{u, n}-y^{u, n}_t)
  +(nt)^{-1/2}e^{\lambda_1t/2}\langle f_1,X_{t}\rangle\nonumber\\
  &=:& J_1^n(t)+J_2^n(t).
\end{eqnarray}
From the proof of \eqref{1.49}, we get that for any $x\in E$,
$$
\P_{\delta_x}J_2^n(t)^2\lesssim n^{-1} (\rho_f^2\phi_1(x)+t^{-1}(a_{t_0}(x)+a_{t_0}(x)^{1/2})).
$$
Thus, there exists $c>0$ such that for any $x\in E$,
\begin{equation}\label{4.3}
 \limsup_{t\to\infty}\P_{\delta_x}J_2^n(t)^2\leq cn^{-1}\phi_1(x).
\end{equation}

Now we consider $J_1^n(t)$.
We define an ${\mathbb R}^2$-valued random variable $U_2(n,t)$ by
\begin{eqnarray*}
U_2(n,t):=\left(e^{\lambda_1 t}\langle \phi_1,X_t\rangle,e^{(\lambda_1/2) t}\sum_{u\in\cL_t}(Y_t^{u, n}-y^{u, n}_t) \right).
\end{eqnarray*}
We claim that
\begin{equation}\label{4.4}
  U_2(n,t)\stackrel{d}{\to}\left(W_\infty, \sqrt{W_\infty}G_2(f)\right), \quad \mbox{ as } t\to\infty.
\end{equation}
Denote the characteristic function of $U_2(n,t)$ under $\P_{\delta_x}$ by
$\kappa_2(\theta_1,\theta_2,n,t)$.
Using an argument similar to that leading to \eqref{2.3}, we get
\begin{eqnarray*}
  \kappa_2(\theta_1,\theta_2,n,t)
  &=& \P_{\delta_x}\left(\exp\{i\theta_1e^{\lambda_1 t}\langle \phi_1, X_t\rangle\}\prod_{u\in\cL_t}h_t^n\left(z_u(t),e^{(\lambda_1/2) t}\theta_2\right)\right),
\end{eqnarray*}
where
$$
h^n_t(x,\theta)=\P_{\delta_x}e^{i\theta(Y^n_t-P_{\delta_x}Y^n_t)}.
$$

Let $t_k,m_k\to\infty$, as $k\to\infty$. Now we consider
\begin{equation}\label{3.16}
  S_k:=e^{\lambda_1t_k/2}\sum_{j=1}^{m_k}(Y_{k,j}-y_{k,j}),
\end{equation}
where $Y_{k,j}$ has the same law as $Y_{t_k}^n$ under $\P_{\delta_{a_{k,j}}}$ and $y_{k,j}=\P_{\delta_{a_{k,j}}}Y_{t_k}^n$. Further, $Y_{k,j},j=1,2,\dots$ are independent. Denote $V_t^n(x):=Var_{\delta_x}Y_t^n$. Suppose the Lindeberg conditions hold:
\begin{description}
  \item{(i)} as $k\to\infty$, $$e^{\lambda_1t_k}\sum_{j=1}^{m_k}\E(Y_{k,j}-y_{k,j})^2=e^{\lambda_1t_k}\sum_{j=1}^{m_k}V_{t^k}^n(a_{k,j})\to\sigma^2;$$
  \item{(ii)} for every $c>0$,
  \begin{eqnarray*}
     && e^{\lambda_1t_k}\sum_{j=1}^{m_k}\E\left(|Y_{k,j}-y_{k,j}|^2,|Y_{k,j}-y_{k,j}|>c e^{-\lambda_1 t_k/2}\right) \\
     &=& e^{\lambda_1t_k}\sum_{j=1}^{m_k}\P_{\delta_{a_{k,j}}}\left(|Y_{t_k}^n-y_{t_k}^n|^2,|Y_{t_k}^n-y_{t_k}^n|>ce^{-\lambda_1 t_k/2}\right)\to 0,\quad k\to\infty.
  \end{eqnarray*}
\end{description}
Then $S_k\stackrel{d}{\to}\mathcal{N}(0,\sigma^2)$ which implies
\begin{equation}\label{3.17}
  \prod_{j=1}^{m_k}h_{t_k}^n(a_{k,j},e^{\lambda_1t_k/2}\theta)\to e^{-\frac{1}{2}\sigma^2\theta^2}, \quad\mbox{as }k\to\infty.
\end{equation}
By Lemma \ref{lem:2.2},
$|V_t^n(x)-\rho_f^2\phi_1(x)|\lesssim ((n-1)t)^{-1}(a_{t_0}(x)^{1/2}+a_{t_0}(x))$
for every $x\in E$.
So by \eqref{1.25}, we get
$$
t^{-1}e^{\lambda_1t}T_t(\sqrt{a_{t_0}}+a_{t_0})(x)\to 0, \mbox{ as } t\to\infty,
$$
which implies
$$
t^{-1}e^{\lambda_1t}\langle \sqrt{a_{t_0}}+a_{t_0},X_t\rangle\to 0, \mbox{ as } t\to\infty,
$$
in probability. Thus,
\begin{equation}\label{3.18}
\lim_{t\to\infty}e^{\lambda_1t}\sum_{u\in \cL_t}V_t^n(z_u(t))=\lim_{t\to\infty}e^{\lambda_1t}\langle\rho_f^2\phi_1,X_t\rangle
=\rho_f^2 W_\infty,\quad\mbox{  in probability}.
\end{equation}
Let
$$
g_n(t,x)=\P_{\delta_x}\left(|Y_t^n-y_t^n|^2,|Y_t^n-y_t^n|>c e^{-\lambda_1 t/2}\right).
$$
We will show that, as $t\to\infty$,
\begin{equation}\label{3.19}
  e^{\lambda_1t}\sum_{u\in \cL_t}\P_{\delta_{z_u(t)}}\left(|Y_t^n-y_t^n|^2,|Y_t^n-y_t^n|>c e^{-\lambda_1 t/2}\right)=e^{\lambda_1t}\langle g_n(t,\cdot),X_t\rangle\to 0,
\end{equation}
in probability.
By Lemma \ref{lem:5.5}, $\lim_{t\to\infty}g_n(t,x)=0$ for every $x\in E$.
Since
$$
g_n(t,x)\leq V_t^n(x)\lesssim \rho_f^2\phi_1(x)+a_{t_0}(x)^{1/2}+a_{t_0}(x)\in L^2(E,\mu),
$$
by the dominated convergence theorem, we have that for any $x\in E$,
$$
\lim_{t\to\infty}\|g_n(t,x)\|_2=0.
$$
By Lemma \ref{lem:rsnew}, we have that for any $x\in E$,
$$
e^{\lambda_1t}\P_{\delta_x}\langle g_n(t,\cdot),X_t\rangle=e^{\lambda_1t}T_t(g_n(t,\cdot))(x)\to 0,\quad \mbox{as}\quad t\to\infty,
$$
which implies \eqref{3.19}.
Thus, for any sequence $s_k\to\infty$, there exists a subsequence $s'_k$ such that,
if we let $t_k=s'_k$, $m_k=|X_{t_k}|$ and $\{a_{k, j}, j=1, \dots, m_k\}=
\{z_u(t_k), u\in {\cal L}_{t_k}\}$, then the Lindeberg conditions hold
$\P_{\delta_x}$-a.s.
Therefore, by \eqref{3.17}, we have
\begin{equation}\label{3.20}
  \lim_{t\to\infty}\prod_{u\in\cL_t}h^n_t(z_u(t),e^{(\lambda_1/2) t}\theta_2)
  =\exp\left\{-\frac{1}{2}\theta_2^2\rho_f^2W_\infty\right\},\quad\mbox{in probability}.
\end{equation}
Hence by the dominated convergence theorem, we get
\begin{equation}\label{3.21}
 \lim_{t\to\infty}\kappa_2(\theta_1,\theta_2,n,t)=\P_{\delta_x}\exp\left\{i\theta_1W_\infty\right\}
  \exp\left\{-\frac{1}{2}\theta_2^2\rho_f^2 W_\infty\right\},
\end{equation}
which implies our claim \eqref{4.4}.
Thus, we easily get that under $\P_{\delta_x}$,
\begin{eqnarray*}
U_3(n,t):=\left(e^{\lambda_1 (nt)}\langle \phi_1, X_{nt}\rangle,J_1^n(t) \right)
\stackrel{d}{\to}\left(W_\infty, \sqrt{\frac{n-1}{n}}\sqrt{W_\infty}G_2(f)\right), \quad \mbox{ as } t\to\infty,
\end{eqnarray*}
where $G_2(f)\sim \mathcal{N}(0,\rho_f^2)$ is independent of $W_\infty$.

Let $\mathcal{L}(nt)$ and $\widetilde{\mathcal{L}}^n(t)$ be the distributions of
$U_1(nt)$ and $U_3(n,t)$
respectively, and let $\mathcal{L}^n$ and $\mathcal{L}$ be the distributions of $(W_\infty, \sqrt{\frac{n-1}{n}}\sqrt{W_\infty}G_2(f))$ and $(W_\infty, \sqrt{W_\infty}G_2(f))$ respectively. Then, using \eqref{5.20}, we have
\begin{eqnarray}\label{4.12}
  \limsup_{t\to\infty}\beta(\mathcal{L}(nt),\mathcal{L})&\leq&
  \limsup_{t\to\infty}[\beta(\mathcal{L}(nt),\widetilde{\mathcal{L}}^n(t))+\beta(\widetilde{\mathcal{L}}^n(t),\mathcal{L}^n)
  +\beta(\mathcal{L}^n,\mathcal{L})]\nonumber\\
 &\leq &\limsup_{t\to\infty}\P_{\delta_x}(J_2^n(t))^2)^{1/2}+0+\beta(\mathcal{L}^n,\mathcal{L}).
\end{eqnarray}
Using this and the definition of $\limsup_{t\to\infty}$, we easily get that
$$
\limsup_{t\to\infty}\beta(\mathcal{L}(t),\mathcal{L})=
\limsup_{t\to\infty}\beta(\mathcal{L}(nt),\mathcal{L})
\le \sqrt{c\phi_1(x)/n}+\beta(\mathcal{L}^n,\mathcal{L}).
$$
Letting $n\to\infty$, we get $ \limsup_{t\to\infty}\beta(\mathcal{L}(t),\mathcal{L})=0$.
The proof is now complete. \hfill$\Box$

\textbf{Proof of Theorem \ref{The:2.3}:}\quad
First note that
\begin{eqnarray*}
  &&t^{-1/2}(\langle \phi_1,X_t\rangle)^{-1/2}\left(\langle f,X_t\rangle-\sum_{\lambda_k<\lambda_1/2}e^{-\lambda_kt}\sum_{j=1}^{n_k}a_j^kH^{k,j}_\infty\right)\\
   &=&  t^{-1/2}(\langle \phi_1,X_t\rangle)^{-1/2}\langle f_{(cl)},X_t\rangle+ t^{-1/2}(\langle \phi_1,X_t\rangle)^{-1/2}\left(\langle f_{(s)}, X_t\rangle-\sum_{\lambda_k<\lambda_1/2}e^{-\lambda_kt}\sum_{j=1}^{n_k}a_j^kH^{k,j}_\infty\right) \\
   &=:& J_1(t)+J_2(t),
\end{eqnarray*}
where $f_{(cl)}=f_{(l)}+f_{(c)}$.
By the definition of $f_{(s)}$, we have $(f_{(s)})_{(c)}=0$. Then using Theorem \ref{The:2.1} for $f_{(s)}$,
we have
\begin{equation}\label{7.2}
  (\langle \phi_1,X_t\rangle)^{-1/2}\left(\langle f_{(s)}, X_t\rangle-\sum_{\lambda_k<\lambda_1/2}e^{-\lambda_kt}\sum_{j=1}^{n_k}a_j^k
  H^{k,j}_\infty\right)\stackrel{d}{\to}G_3(f_{(s)}).
\end{equation}
 Thus
 \begin{equation}\label{7.3}
 J_2(t)\stackrel{d}{\to}0,\quad t\to\infty.
 \end{equation}
Since $\lambda_1=2\lambda_{\gamma(f_{(cl)})}$, so using Theorem \ref{The:1.4} for $f_{(cl)}$, we have
\begin{equation}\label{7.1}
(e^{\lambda_1 t}\langle \phi_1,X_t\rangle, J_1(t)))\stackrel{d}{\to}(W^*, G_2(f_{(cl)})),
\end{equation}
where $G_2(f_{(cl)})\sim\mathcal{N}(0,\rho_{f_{(cl)}}^2)$.
By the definition of $\rho^2_f$ given by \eqref{e:rho}, we have $\rho_{f_{(cl)}}^2=\rho_{f_{(c)}}^2$.
Combining \eqref{7.3} and \eqref{7.1}, we arrive at the conclusion of Theorem \ref{The:2.3}.\hfill$\Box$

\vspace{.1in}
\begin{singlespace}

\end{singlespace}

\vskip 0.2truein
\vskip 0.2truein

\noindent{\bf Yan-Xia Ren:} LMAM School of Mathematical Sciences \& Center for
Statistical Science, Peking
University,  Beijing, 100871, P.R. China. Email: {\texttt
yxren@math.pku.edu.cn}

\smallskip
\noindent {\bf Renming Song:} Department of Mathematics,
University of Illinois,
Urbana, IL 61801, U.S.A.
Email: {\texttt rsong@math.uiuc.edu}

\smallskip

\noindent{\bf Rui Zhang:} LMAM School of Mathematical Sciences, Peking
University,  Beijing, 100871, P.R. China. Email: {\texttt
ruizhang8197@gmail.com}

\end{doublespace}
\end{document}